\documentclass[11pt,letterpaper,headinclude,footinclude,fleqn,reqno]{amsart}                 
\usepackage[T1]{fontenc}                   
\usepackage[utf8]{inputenc}                 
\usepackage[english]{babel}       
\usepackage{graphicx}                      %
\usepackage[font=small]{quoting}            %
\usepackage{caption}  
\usepackage{amsmath}

\usepackage[top=.8in,bottom=1.2in,left=1.5in,right=1.5in]{geometry}
\usepackage{verbatim}
\usepackage{color}
\usepackage{picture}
\usepackage{graphicx}
\usepackage{graphics}
\usepackage{comment}
\usepackage{hyperref}
\hypersetup{colorlinks,linkcolor={blue},citecolor={blue},urlcolor={red}}

\def\n{\nabla}
\def\<{\langle}
\def\>{\rangle}

\def\O{\Omega}

\def\S{\Sigma}
\def\n{\nabla}

\def\n{\nabla}

\def\O{\Omega}
\def\p{\partial}

\def\s{\sigma}

\def\arr{\rightarrow}

\def\n{\nabla}
\def\<{\langle}
\def\>{\rangle}

\def\n{\nabla}

\def\RR{\mathbb{R}}

\def\SS{\mathbb{S}}
\def\HH{\mathbb{H}}
\def\hh{\mathbb{H}}
\def\BB{\mathbb{B}}
\def\RR{\mathbb{R}}
\def\rr{\mathbb{R}}
\def\Si{\mathcal{S}}

\def\O{\Omega}
\def\p{\partial}

\def\s{\sigma}

\def\arr{\rightarrow}

\newtheorem{theorem}{Theorem}[section]

\newtheorem{proposition}[theorem]{Proposition}
\newtheorem{corollary}[theorem]{Corollary}

\theoremstyle{definition}
\newtheorem{definition}[theorem]{Definition}

\usepackage{graphicx}

\numberwithin{equation}{section}

\makeatletter
\@namedef{subjclassname@2020}{\textup{2020} Mathematics Subject Classification}
\makeatother

\begin{document}
	\setlength{\baselineskip}{1.2\baselineskip}

	\title[A locally constrained mean curvature type flow]
	{ A constrained mean curvature type flow for capillary boundary  hypersurfaces  in space forms }

	\author{Xinqun Mei}
	\address{School of Mathematical Sciences, University of Science and Technology of China, Hefei, 230026, Anhui Province, P.R. China}

    \address{Mathematisches Institut, Albert-Ludwigs-Universit\"{a}t Freiburg, Freiburg im Breisgau, 79104, Germany}

	\email{qunmx@mail.ustc.edu.cn}

\author{Liangjun Weng} 
		
\address{Dipartimento Di Matematica, Universita Di Roma Tor Vergata, Via della Ricerca Scientifica, I-00133 Roma, Italy}
	
\email{ljweng08@mail.ustc.edu.cn}
		
	\subjclass[2020]{	Primary 53C44, Secondary 35K93.}
	\keywords{Mean curvature type flow,  space forms, capillary boundary, the capillary isoperimetric inequality. }
	
	\maketitle
	\begin{abstract}
		In this paper, we introduce a new constrained mean curvature type flow for capillary boundary  hypersurfaces   in space forms. %
		We show the flow exists for all time and converges globally to a spherical cap. Moreover, the flow preserves the volume of the bounded  domain enclosed by the hypersurface and decreases the total energy. As a by-product, we give a flow proof of the capillary isoperimetric inequality for the  star-shaped capillary boundary hypersurfaces in   space forms.

	\end{abstract}
	
	\bigskip

	\section{Introduction}
	Let $M^{n+1}(K)$ be a complete simply-connected Riemann manifold with constant sectional curvature $K$. Up to  homoteties, we may assume that $K=-1, 0, +1$.  The case $K=0$, $M^{n+1}(K)$ is just the Euclidean space $\RR^{n+1}$. If  $K=-1$, $M^{n+1}(K)$ is the hyperbolic space $\HH^{n+1}$, and we use the  Poincar\'e ball model for   $\hh^{n+1}$, which  is given by $\left(\BB^{n+1}, \bar g_{\HH} \right)$, where
	\begin{eqnarray*}
		\BB^{n+1}=\{x\in\rr^{n+1}: |x|<1\}, \quad \bar g_{\HH}= e^{2u}\delta_{\BB^{n+1}}:=\frac{4}{(1-|x|^2)^2}|dx|^2.
	\end{eqnarray*}
	Let $B_R^{\mathbb{H}}$ be a geodesic ball in $\mathbb{H}^{n+1}$ with the hyperbolic radius $R\in (0,+\infty)$, by using the isometry of $\mathbb{H}^{n+1}$, one can view $B_R^{\mathbb{H}}$ as a Euclidean ball $B_{r_0}\subset \mathbb{B}^{n+1}$ of  radius 
	$$ { r_0:=\sqrt{\frac{\cosh  R-1}{ \cosh R+1}}}\in (0,1),$$ in $\mathbb{R}^{n+1}$ with the hyperbolic metric $\bar g_{\HH}$. 
	If $K=1$, $M^{n+1}(K)$ is the spherical space  form $\SS^{n+1}$, and we use the model
    \begin{eqnarray*}
    	(\RR^{n+1}, \bar{g}_{\SS}), \quad  \bar{g}_{\SS}= e^{2v}\delta_{\RR^{n+1}}:=\frac{4}{(1+|x|^2)^2}|dx|^2.
    \end{eqnarray*}
    to represent $\SS^{n+1}\setminus \{\Si\}$, the unit ball without the south pole. Let $B_{R}^{\SS}$ be a geodesic ball in $\SS^{n+1}$ with radius $R\in (0, \pi)$ center at the north pole. Analog to hyperbolic space, one can view $B_{R}^{\SS}$ as a Euclidean ball $B_{r_{0}}\subset \RR^{n+1}$ of radius 
    $$r_{0}:=\sqrt{\frac{1-\cos R}{1+\cos R}}\in (0, \infty).$$
   with the spherical metric $\bar{g}_{\SS}$. 
   For brevity, and without causing ambiguity, we uniformly denote $\bar{g}$  for  the metric $\bar{g}_{\HH}$ or $\bar{g}_{\SS}$ and $B_{R}$ for the geodesic ball $B_{R}^{\HH}$ in $\HH^{n+1}$ or $B_{R}^{\SS}$ in $\SS^{n+1}$ with the geodesic radius $R$ in the rest of this paper.
	
	For a constant unit vector $a\in \SS^{n}\subset \RR^{n+1}$, let $X_a$ be the vector field in $M^{n+1}(K)$ $(K=-1, +1)$, given by
	\begin{eqnarray*}\label{conf killing}
		X_a:=\frac{2}{1+Kr_0^2}\left[\langle x, a\rangle x-\frac{1}{2}(|x|^2+r_0^2)a\right],
	\end{eqnarray*}
	$\< ,\>$ denotes the Euclidean metric. It has been observed by Wang-Xia \cite[Proposition 4.1]{WX19} that $X_a$ is a conformal Killing vector field  such that 
	$$L_{X_a}\bar g =V_a \bar g,$$
	where 
	\begin{eqnarray*}\label{v e}
		V_a:=\frac{2\langle x, a\rangle }{1+K|x|^2}.
	\end{eqnarray*}
    Besides, $X_a$ satisfies $\bar g(X_a,\bar N)=0$ on $\p B_R$, where   $\bar N$ is the  unit normal vector along $\p B_R$. From those properties of $X_{a}$, Wang-Xia \cite[Proposition 4.4]{WX19} obtain a new Minkowski formula in space forms. That is, for a hypersurface $\S\subset \bar {B}_R$ with boundary $\p\S\subset \p B_R$ such that $\p \S$ intersects $\p B_R$   at a constant contact angle $\theta\in (0,\pi)$, in $\HH^{n+1}$, it holds 
	\begin{eqnarray}\label{minkowski hyperbolic}
		\int_{\Sigma}n\left( V_{a}+ \sinh R\cos \theta\bar{g}(Y_{a},\nu)\right)dA=\int_{\Sigma}H\bar{g}(X_{a},\nu)dA,
	\end{eqnarray}
    and in $\SS^{n+1}$ respectively,
    \begin{eqnarray}\label{minkowski sphere}
    		\int_{\Sigma}n\left( V_{a}+ \sin R\cos \theta\bar{g}(Y_{a},\nu)\right)dA=\int_{\Sigma}H\bar{g}(X_{a},\nu)dA,
    \end{eqnarray} 
	where $H\nu$ is mean curvature vector of $\Sigma$ and $Y_a:=\frac{1}{2}(1-K|x^{2}|)a+K\<x,a\>x$.

	\
In this paper, 
we  consider a new type  mean curvature flow for capillary boundary hypersurfaces (see Section \ref{sec2.1} for the definition) supported in the geodesic ball in space forms $M^{n+1}(K)$ for $K=-1$ and $K=1$, while the case $K=0$ was  studied previously in \cite{WW20,WX20} respectively. 	Such kind of locally constrained curvature type flow was first used
by Guan-Li \cite{GL15} for closed hypersurfaces in space forms, which was motivated by the Minkowski formula. See also \cite{GL18,GLW,HL,HLW,LS21,SX19} and references therein for various general setting, which include the hyperbolic space $\HH^{n+1}$ and spherical space $\SS^{n+1}$. 

	Recently, it attracts high interest to study the flow of hypersurface with non-empty boundary, especially due to their close connection with geometric inequalities. For example, some new class of constrained curvature flows was studied for hypersurface with free boundary in  Euclidean space by \cite{SWX,WX20}, and  \cite{MWW,WW20,WeX21} for capillary boundary. Subsequently, a class of new Alexandrov-Fenchel inequalities was obtained after establishing the long-time existence and convergence of those flows. One can refer to \cite{LS16,LS17} for the studying of the inverse mean curvature flow with free boundary in a Euclidean ball and its application in geometric inequality. 

Therefore, besides the Euclidean space, it is natural to ask the same question in general ambient space, say space forms.		
To be more precise, let $\Sigma_{t}$ be a family of hypersurfaces with boundary in $\bar {B}_{R}$ given by a family of isometric embeddings $x(\cdot, t): M \rightarrow  \bar {B}_{R}$ from a compact $n$-dimensional manifold $M$ with the boundary $\partial M$ ($n\geq 2$) such that 
	\begin{eqnarray}
		{\rm 	int}(\Sigma_{t})=x\left({\rm int}(M), t\right)\subset B_{R},\quad \partial\Sigma_{t}=x(\partial M,t)\subset \partial B_{R}.\notag
	\end{eqnarray}
And $x(\cdot, t)$ satisfy
	\begin{equation}\label{flow with capillary in hyperbolic}
		\left\{ \begin{array}{lll}
			&(\partial_t x)^\perp(\cdot, t)  =F(\cdot, t) \nu(\cdot, t), \quad &
			\hbox{ in }M\times[0,T),\\
			&\bar g( \nu(\cdot, t),\bar N\circ x(\cdot,t)) = -\cos\theta 
			\quad & \hbox{ on }\partial M \times [0,T),\\
			& x(\cdot,0)  = x_0(\cdot) \quad & \text{ in }   M,
		\end{array}\right.
	\end{equation}
	
	In the case of $\HH^{n+1}$, we choose the speed function in \eqref{flow with capillary in hyperbolic} as
	\begin{eqnarray}\label{speed f hyperbolic}
F:=nV_{a}+n\sinh R\cos \theta \bar{g}(Y_{a},\nu)-H\bar{g}(X_{a}, \nu).
	\end{eqnarray}
 In the case of $\SS^{n+1}$, we choose the speed function   in  \eqref{flow with capillary in hyperbolic} as 
 \begin{eqnarray}\label{speed f sphere}
F:=nV_{a}+n\sin R\cos \theta \bar{g}(Y_{a},\nu)-H\bar{g}(X_{a}, \nu).
 \end{eqnarray}

Such kind of flow is motivated by the Minkowski formula \eqref{minkowski hyperbolic} and \eqref{minkowski sphere} respectively. And the speed function $F$ is chosen to be \eqref{speed f hyperbolic} or \eqref{speed f sphere} to ensure that the enclosed volume of $\S_t\subset \bar B_R$ is preserved along the flow \eqref{flow with capillary in hyperbolic},  
while the total energy functional \eqref{wetting energy}  is monotone decreasing.  We will discuss this later.

	\begin{definition}\label{star-shaped}
		We say   $\Sigma\subset  \bar {B}_R\subset M^{n+1}(K)$ is star-shaped with respect to $a$ 
		if $\bar g( X_a,\nu) >0$ along $\S$.
	\end{definition}
	Our main result in this paper is the following theorem.
	\begin{theorem}\label{thm1.1}
		If the initial hypersurface $\S_0$ is  star-shaped capillary boundary hypersurface  in space forms $M^{n+1}(K)$ 
		and the contact angle satisfies $|\cos\theta|<\frac{3n+1}{5n-1}$, then the flow \eqref{flow with capillary in hyperbolic} 
		exists for all time. Moreover, $x(\cdot, t)$ converges to a spherical cap in the $C^{\infty}$ topology as $t\rightarrow \infty$, whose enclosed domain has the same volume as the one enclosed by $\Sigma_{0}$.
	\end{theorem}
	The family of the spherical cap is given by geodesic ball of radius $r$ and totally geodesic ball in the hyperbolic space or in the spherical space, which can be viewed as a Euclidean set
	\begin{eqnarray}\label{static model}
		C_{ \theta,r}(a):=\left\{x\in B_R: \left|x- \sqrt{r^2+2rr_0\cos\theta+r_0^2} a\right|= r \right\}, 
	\end{eqnarray}and 
	\begin{eqnarray}
		C_{ \theta,\infty}(a):=\{x\in  B_R  :  \<x, a\>=\cos\theta\},
	\end{eqnarray}endowed with the  metric $\bar g$.

	In particular, when $K=0$ or  $\theta=\frac{\pi}{2}$ and $K=-1$, Theorem \ref{thm1.1} was proved in \cite{WW20} and \cite{QWX} respectively.  And we have a technique restriction on the range of contact angle $|\cos\theta|<\frac{3n-1}{5n-1}$ as in \cite{WW20}, which is crucial for us to obtain uniform gradient estimate in space forms, see Section \ref{sec3}. However,  we expect the result holds true for the whole range of $\theta\in (0,\pi)$ in space forms.
	
	It is worth noting that the isoperimetric inequality for hypersurfaces with non-empty boundary in a Euclidean ball were studied in \cite{BS,BM}, and they proved that among the hypersurface with fixed volume of enclosed domain  $\O$, the spherical cap are the minimizer of area   functional $\S$. Instead of just considering the area functional, it is interesting to consider the  total energy  functional, that is
	\begin{eqnarray}\label{wetting energy}  
		E(\Sigma):=	{\rm Area}(\S)-\cos \theta {\rm Area}(T)
	\end{eqnarray} 
	for a  constant contact angle $\theta \in (0,\pi)$. The second term $T:=\p\O\setminus \S$ is known as the wetting part of $\p\O$ in the theory of capillarity, see \cite{Finn} for example.	
	
	Further, by combining Theorem \ref{thm1.1} with our high order Minkowski formulas  \eqref{minkowski k}   and \eqref{minkowski formula sphere k}
	for $k=2$, we give a flow proof for the capillary isoperimetric inequality  in space forms for \eqref{wetting energy}, which can also be viewed as the hyperbolic  and sphere counterpart in \cite[Theorem 1.1]{WW20} or \cite[Chapter 19]{Maggi}.  
	\begin{corollary}\label{cor1.3}
		Among the star-shaped capillary boundary hypersurfaces with fixed volume of enclosed domain  in a geodesic ball $ \bar{B}_{R}\subset M^{n+1}(K)$ for $K=\pm 1$, the spherical caps are the only minimizers of the total energy \eqref{wetting energy},
		provided that the contact angle $\theta $ satisfies $|\cos\theta|<\frac{3n+1}{5n-1}$.
	\end{corollary}

	\begin{proof}
For $K=-1$, along the flow \eqref{flow with capillary in hyperbolic} with $F$ being  \eqref{speed f hyperbolic}, using \eqref{2.4} and \eqref{minkowski hyperbolic},  we know
		\begin{eqnarray*}
			\frac{d}{dt}{\rm Vol}({\O_t})&=&\int_{\Sigma_{t}}\left[nV_{a}+n\sinh R\cos\theta \bar g(Y_{a}, \nu)-H\bar
			{g}(X_{a},\nu) \right]dA_{t}\\&=&0.
		\end{eqnarray*} 
		Note that $\p_t x\big|_{\p M}\in T (\p B_R^{\HH})$, combining with \eqref{2.5} and \eqref{minkowski k} for $k=2$, it follows   
		\begin{eqnarray*}
			\frac{d}{dt} E(\S_t)&=&\int_{\Sigma_{t}}H\left[nV_{a}+n\sinh R\cos \theta \bar{g}(Y_{a},\nu)-H\bar{g}(X_{a}, \nu)\right]dA_{t}\\
			&=&\int_{\Sigma_{t}}(\frac{2n}{n-1}\sigma_{2}-H^{2})\bar{g}(X_{a},\nu)dA_{t}\\
			&=&-\frac{1}{n-1}\int_{\Sigma_{t}}\sum\limits_{1\leq i< j\leq n}(\kappa_{i}-\kappa_{j})^{2} \bar{g}(X_{a},\nu)dA_{t}\leq 0.
		\end{eqnarray*} 
For $K=1$, the proof is similar to above.
	Hence Corollary \ref{cor1.3} follows directly from Theorem \ref{thm1.1}.
	\end{proof}
	
	\vspace{.2cm}
	\textbf{This article is structured as follows.} In Section \ref{sec2}, we give some preliminaries for  capillary boundary hypersurfaces supported in a geodesic ball of space forms, and prove the high-order Minkowksi type formula in the hyperbolic space $\HH^{n+1}$ and in the spherical space $\SS^{n+1}$ respectively. Then we  convert the flow \eqref{flow with capillary in hyperbolic} 
	to a scalar parobolic equation  on semi-sphere with the help of a conformal transformation. The last Section \ref{sec3} is devoted to obtain uniform a priori estimates and prove  Theorem \ref{thm1.1}.
	
	\section{Preliminaries}\label{sec2}

	In this section, in the first part, we collect some basic facts about capillary boundary hypersurfaces supported in a geodesic ball of space forms, and we establish Minkowski formula which can be viewed as a high order counterpart of \eqref{minkowski hyperbolic} and \eqref{minkowski sphere}. In the second part, we reduce \eqref{flow with capillary in hyperbolic} 
	to a scalar  flow, provided that the evolving hypersurface is star-shaped in the sense of Definition \ref{star-shaped}.

	\subsection{Capillary boundary hypersurfaces in space forms}\label{sec2.1}
	
Let $x:M\to \bar B_R$ be an isometric embedding of an orientable $n$-dimensional compact manifold $M$ with smooth boundary $\p M$, denote $\S:=x(M)$ and $\p \S:=x(\p M)$ which satisfying  
	${\rm int}(\Sigma)\subset {\rm int}(B_{R})$ and $\partial\Sigma\subset \partial B_{R}$. If there is no confusion, we will identity $M$ with $\Sigma$ and $\partial M$ with $\partial \Sigma$. We denote by $\bar{\nabla}, D$ the Levi-Civita connection of $(B_{R}, \bar g)$ and $(\Sigma, g)$ respectively, where $g$ is the induced metric from embedding $x$.

	$\Sigma$ divides the  ball $B_{R}$ into two parts, we  denote one part by $\Omega$ and $\nu$
	be the unit outward normal vector of $\Sigma$ w.r.t $\Omega$ and $T:=\p\O\cap \p B_R$ be the wetting part of $\p\O$. Let $\mu$ be the unit outward conormal vector field along $\partial\Sigma\subset \S$ and $\bar\nu$ be the unit normal to $\partial\Sigma\subset \partial B_{R}$ such that $\{\nu, \mu\}$ and $\{ \bar{\nu}, \bar{N}\}$ have the same orientation in the normal bundle of $\partial\Sigma\subset {B}_{R}$. Denote by $h$ and $\sigma_{k} $ the second fundamental form and $k$-th mean curvature of the immersion $x$ respectively, precisely, $h(X,Y):=\bar{g}(\bar{\nabla}_{X}\nu, Y)$ and $\sigma_{k}:=\sigma_{k}( \kappa)$, where $X,Y\in T\S$ and  $\kappa:=(\kappa_1,\cdots,\kappa_n)\in\RR^n$ are the eigenvalues of Weingarten matrix $  (h^j_i)$. 
	We define the {\it contact angle} $\theta\in (0,\pi)$ between  the hypersurface $\Sigma$ and the geodesic ball $B_R$ by
	$$\bar g( \nu, \bar N)=\cos (\pi-\theta),\quad \text{ on } \p\S.$$ (See for example Figure 1 in the case $\HH^{n+1}$.) And we call such hypersurface $\S$ as the {\it capillary boundary hypersurface}. In particular, when $\theta=\frac{\pi}{2}$, it is known as {\it free boundary hypersuface}. Moreover, it follows that
	\begin{eqnarray*}	\mu&=&\sin\theta\bar{N}+\cos\theta\bar{\nu},\\		 \nu&=&-\cos\theta \bar{N}+\sin\theta\bar{\nu}.\end{eqnarray*}

	\begin{center}
		\begin{figure}[h] \includegraphics[width=0.45\linewidth]{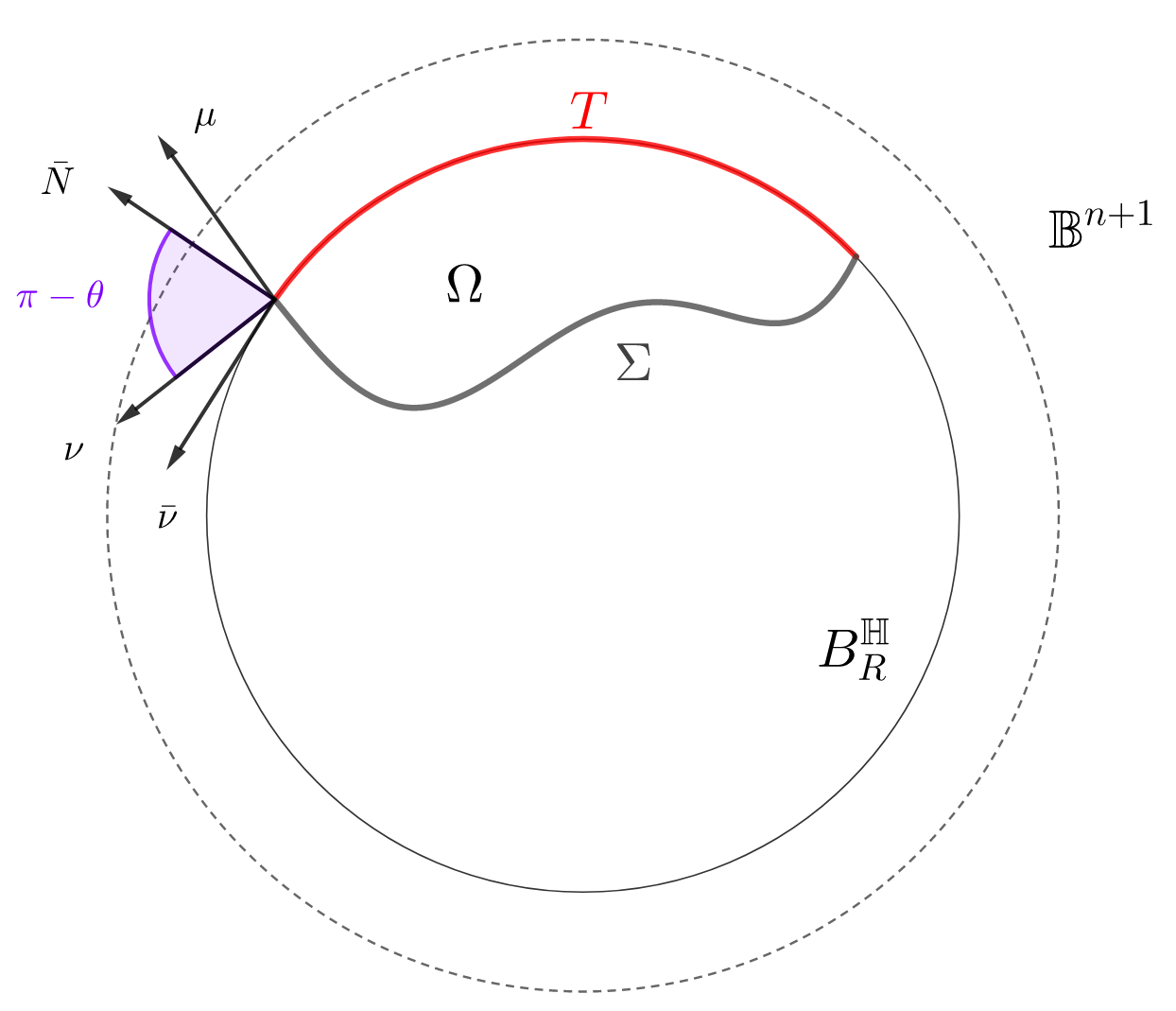} \label{fig1}  \caption*{Figure 1. A capillary boundary hypersurface $\S$   in   hyperbolic space}   
		\end{figure}
	\end{center}

	The so-called wetting energy $W(\S)$ is just the area of 
	the region $T$, which is   bounded by $\partial \Sigma$ on $\partial B_R$. And the total energy functional is defined as
	\begin{eqnarray}\label{total energy}
		E(\Sigma):= \mbox{Area}(\Sigma)-\cos\theta\,  W(\S).
	\end{eqnarray}
	In order to show the monotonicity property of total energy functional $E$ along our flow \eqref{flow with capillary in hyperbolic}. 
	Let us consider an  admissible variation of $\S:=x(M)$, given by $x_{t}: M \times (-\varepsilon, \varepsilon)\rightarrow \bar{B}_{R}$ 
	satisfying that $x_{t}(\cdot):=x(\cdot, t): M \rightarrow \bar{B}_{R}$ is an immersion with  
	\begin{eqnarray*}
		{\rm int}(\Sigma_{t}) \subset  B_{R},\quad  {\rm and  }\quad  \partial\Sigma_{t} \subset \partial B_{R},
	\end{eqnarray*}
	where $\S_t:=x(M,t)$. Denote $\O_t$ the enclosed domain by $\S_t$ in $\bar B_R$.
	
	Let $Y:=\frac{\p}{\p t} x(\cdot,t)\big|_{t=0}$ be the associated variational vector field of $x_t$, then the first variational formula of ${\rm Vol}(\O_t)$ and $E(\S_t)$ are known as (cf. \cite[Section 4]{RR})
	\begin{eqnarray}\label{2.4}
		\frac{d}{dt}\Big|_{t=0}{\rm Vol}(\O_t)=\int_{M}\bar g(Y, \nu)d A,
	\end{eqnarray}	
	and
	\begin{eqnarray}\label{2.5}
		\frac{d}{dt}\Big|_{t=0}E(\S_t)=\int_{M}H \bar g(Y,\nu)dA+\int_{\partial M}\bar g(Y, \mu-\cos\theta \bar\nu)ds,
	\end{eqnarray}
	where $dA$ and $ds$ are the area element of $M$ and $\partial M$ respectively.
	
	As the Minkowski formula plays an important role for the closed hypersurface  in space forms, cf.\cite{GL15,GL18,GLW} etc. In the following, we firstly establish  a new Minkowski formula for the capillary boundary hypersurfaces supported in a geodesic ball in space forms $M^{n+1}(K)$ for $K=\pm 1$. While the case $K=0$ was shown recently by Weng-Xia in \cite[Proposition 2.8]{WeX21}. It may also have independent interest for $K=\pm 1$, since we only use the special case $k=2$ of \eqref{minkowski k} and \eqref{minkowski formula sphere k}  to show the monotonicity of the total energy functional \eqref{total energy} in this paper, which has been indicated in the proof of Corollary \ref{cor1.3}.
	\begin{proposition}\label{pro2.1}  
		Let $x:M\rightarrow B_{R}^{\mathbb{H}}$ be an isometric immersion of $\Sigma:=x(M)$ into the hyperbolic ball $B_{R}^{\mathbb{H}}$, whose boundary $\partial \Sigma$ intersects $\partial B_{R}^{\mathbb{H}}$ at a constant angle $\theta\in (0,\pi)$, then it holds
		\begin{eqnarray}\label{minkowski k}
			(n-k+1)\int_{\Sigma}\big[\sigma_{k-1} V_{a}+\sinh R\cos\theta \sigma_{k-1}\bar g(Y_{a},\nu)\big]dA=k\int_{\Sigma}\sigma_{k}\bar g(X_{a},\nu)dA. 
		\end{eqnarray}
	\end{proposition}
	In particular, when $\theta=\frac{\pi}{2}$, formula \eqref{minkowski k} was proved by Wang-Xia \cite[Proposition 5.1]{WX19}. For completeness, we contain a proof here for general $\theta$.
	\begin{proof}
		Let $\{e_{i}\}_{i=1}^{n}$ be an othonormal frame on $\Sigma$, by  \cite[equation (4.7)]{WX19}, we have
		\begin{eqnarray}\label{2.16}
			\frac{1}{2}(D_{i}(X_{a}^{T})_{j}+D_{j}(X_{a}^{T})_{i})=V_{a}\bar g_{ij}-h_{ij}\bar g(X_{a},\nu),\notag
		\end{eqnarray}
		where $X_{a}^{T}$ is the tangential projection of $X_{a}$ on $\Sigma$. Set
		\begin{eqnarray}
			Z_{a}:=\bar g(\nu, e^{-u}a)x-\bar g(x,\nu)(e^{-u}a).\notag
		\end{eqnarray}
		It is known that $x$ is a conformal Killing vector \cite[Lemma~2.1]{GL15}, satisfying
		\begin{eqnarray}\label{deriv 1}
			\bar{\n}x=V_{0}\bar{g},
		\end{eqnarray} 
		where $V_{0}:=\cosh R =\frac{1+|x|^{2}}{1-|x|^{2}}$ and $R:=R(x)$ is the hyperbolic distance to the origin.  By \cite[Proposition~4.3]{WX19},      \begin{eqnarray}\label{deriv 2}
			\bar{\n}_{Z}(e^{-u}a)=e^{-u}\left[\bar{g}(x, e^{-u}a)Z-\bar{g}(Z, e^{-u}a)x\right],
		\end{eqnarray}
		for any $Z\in T\S$. Using \eqref{deriv 1} and \eqref{deriv 2}, it yields
		\begin{eqnarray*}
			D_{i}(\bar g(\nu, e^{-u}a)x^{T})_{j}&=&
			e^{-u}h_{ik}\bar g(e_{k},a)\bar g(x,e_{j})
			-e^{-u}\bar g(e_{i}, e^{-u}a)\bar g(x, e_{j})\bar g(\nu, x)\\
			&&+\bar{g}(\nu, e^{-u}a)\left[V_{0}\bar g_{ij}-h_{ij}\bar g(x,\nu)\right],
		\end{eqnarray*}
		and 
		\begin{eqnarray*}
			D_{i}(\bar g(x,\nu)(e^{-u}a)^{T})_{j}&=&
			e^{-u}h_{ik}\bar g(x,e_{k})\bar g(a, e_{j})+\bar g(x,\nu)\big[e^{-u}\big(\bar g(x,e^{-u}a)\bar g_{ij}\\
			&&-\bar g(e_{i}, e^{-u}a)\bar g(x, e_{j})\big)-h_{ij}\bar g(e^{-u}a,\nu)\big].\notag
		\end{eqnarray*}
		On the other hand, we have
		\begin{eqnarray}
			\bar g(Z_{a}^{T},\mu)|_{\partial \Sigma}=r_{0}\bar g(\bar \nu, a),\notag
		\end{eqnarray}
		and 
		\begin{eqnarray}
			\bar g(X_{a}^{T},\mu)|_{\partial\Sigma}=-\frac{2r_{0}^{2}}{1-r_{0}^{2}}\cos\theta \bar g(a,\bar \nu),\notag
		\end{eqnarray}
		combining with $\sinh R=\frac{2r_{0}}{1-r_{0}^{2}}$, it follows
		\begin{eqnarray}\label{2.5.}
			\bar g(X_{a}^{T}+\sinh R\cos\theta Z_{a}^{T}, \mu)|_{\partial \Sigma}=0.
		\end{eqnarray}
		Denote $\sigma_{k-1}^{ij}:=\frac{\partial \sigma_{k}}{\partial h_{j}^{i}}$ be the $k$-th Newton transformation,  note that
		\begin{eqnarray*}
			V_{0}\bar{g}(\nu, e^{-u}a)-e^{-u}\bar{g}(x,\nu)\bar{g}(x, e^{-u}a)&=&
			\bar{g}\left(\nu, \frac{1}{2}(|x|^{2}+1)a-\<x, a\>x\right)\\
			&=&\bar{g}(\nu,Y_{a} ),
		\end{eqnarray*} we have
		\begin{eqnarray*}
			&&\sigma_{k-1}^{ij}D_{i}(X_{a}^{T}+\sinh R\cos\theta Z_{a}^{T})_{j}\\
			&=&(n-k+1)\sigma_{k-1}V_{a}-k\sigma_{k}\bar g(X_{a},\nu)+\sinh R\cos\theta\bar{g}(\nu, e^{-u}a)V_{0}\sigma_{k-1}^{ij}\bar{g}_{ij}\\
			&&-\sinh R\cos\theta e^{-u}\bar{g}(x,\nu)\bar{g}(x, e^{-u}a)\sigma_{k-1}^{ij}\bar{g}_{ij}\\
			&=&(n-k+1)\sigma_{k-1}V_{a}-k\sigma_{k}\bar g(X_{a},\nu)+(n-k+1)\sinh R\cos\theta \sigma_{k-1}
			\big[V_{0}\bar{g}(\nu, e^{-u}a)\\
			&&- e^{-u}\bar{g}(x,\nu)\bar{g}(x, e^{-u}a)\big]\\
			&=&(n-k+1)\sigma_{k-1}V_{a}-k\sigma_{k}\bar g(X_{a},\nu)+(n-k+1)\sinh R\cos\theta \sigma_{k-1}\bar g(Y_{a},\nu),
		\end{eqnarray*} 
		integrating above identity on $\Sigma$ and by the divergence theorem, we have
		\begin{eqnarray*}
			&&\int_{\Sigma}\left[(n-k+1)\sigma_{k-1}V_{a}+(n-k+1)\sinh R\cos\theta\sigma_{k-1}\bar{g}(Y_{a},\nu)-k\sigma_{k}\bar{g}(X_{a}, \nu) \right]dA\\	
			&=&\int_{\Sigma}\sigma_{k-1}^{ij}D_{i}(X_{a}^{T}+\sinh R\cos\theta Z_{a}^{T})_{j}dA=
			\int_{\Sigma}D_{i}\left(\sigma_{k-1}^{ij}\bar g(X_{a}^{T}+\sinh R\cos\theta Z_{a}^{T},e_{j})\right)dA\\
			&=&\int_{\partial\Sigma}\sigma_{k-1}^{ij}\bar g(X_{a}^{T}+\sinh R\cos\theta Z_{a}^{T}, e_{j})\overline{g}(\mu, e_{i})ds=\int_{\partial \Sigma}\sigma_{k-1}^{\mu\mu}\bar g(X_{a}^{T}+\sinh R\cos\theta Z_{a}^{T}, \mu)ds =0,
		\end{eqnarray*}
		where the last line follows from that $\mu$ is the principal direction of $\p\S\subset \S$ (cf. 
		\cite[Propostion 2.1]{WX19})  and \eqref{2.5.}.  Thus we complete the proof of Proposition \ref{pro2.1}.
	\end{proof}	
Except the case $\HH^{n+1}$,  we obtain the corresponding Minkowski formula for capillary boundary hypersurfaces supported in a geodesic ball of spherical space $\SS^{n+1}$.
    \begin{proposition}\label{pro minkowski formula}
    	Let $x: M\arr B_{R}^{\SS}$ be an isometric immersion of $\Sigma:=x(M)$ into a ball $B_{R}^{\SS}$ in $\SS^{n+1}$, whose boundary $\partial \Sigma$ intersects $\partial B_{R}^{\SS}$ at a constant angle $\theta\in(0, \pi)$, then it holds 
    	\begin{eqnarray}\label{minkowski formula sphere k}
    		\int_{\Sigma}\left[(n-k+1) \sigma_{k-1}V_{a}+(n-k+1)\sin R\cos\theta\sigma_{k-1}\bar{g}(Y_{a},\nu)\right]dA=k\int_{\Sigma}\sigma_{k}\bar{g}(X_{a}, \nu)dA.
    	\end{eqnarray}
    \end{proposition}The proof is similar to Proposition \ref{pro2.1}, we leave it to the interested readers. 

	\subsection{A scalar equation }
	In this section, we firstly reduce the flow \eqref{flow with capillary in hyperbolic} to a scalar equation on $\bar{\mathbb{S}}^{n}_{+}$, if the evolving hypersurfaces are
	star-shaped in the sense of Definition \ref{star-shaped}.
	With loss of generality, we assume $a:=-e_{n+1}=(0,\cdots, 0, -1)$ from now on, and the half-space
	\begin{eqnarray}
		\mathbb{R}^{n+1}_{+}:=\{x|x=(x_{1},\cdots, x_{n+1})\in \mathbb{R}^{n+1}: x_{n+1}>0\}.\notag
	\end{eqnarray}
	We use the polar coordinate $(\rho,\beta,\theta)\in [0,+\infty)\times [0, \frac{\pi}{2}]\times \mathbb{S}^{n-1}$ in $\mathbb{R}^{n+1}_{+}$, then the standard Euclidean metric in $\mathbb{R}^{n+1}_{+}$ has the form
	\begin{eqnarray}
		|dz|^{2}=d\rho^{2}+\rho^{2}g_{\overline{\mathbb{S}}^{n}_{+}}=d\rho^{2}+\rho^{2}(d\beta^{2}+\sin^{2}\beta g_{\mathbb{S}^{n-1}}).\notag
	\end{eqnarray}
	
	Let $f$  the conformal diffeomorphism from the unit ball to the half space  (cf. \cite[Section 3.2]{WW20} or \cite{QWX,WX19}) given by
	\begin{eqnarray}
		f:&\bar {\mathbb{R}}^{n+1}_{+} &\rightarrow \bar {\mathbb{B}}^{n+1},\notag\\
		&(z^{'}, z_{n+1}) &\mapsto (\frac{2z^{'}}{|z^{'}|^{2}+(1+z_{n+1})^{2}}, \frac{|z|^{2}-1}{|z^{'}|^{2}+(1+z_{n+1})^{2}}):=(y^{'}, y_{n+1}).\notag
	\end{eqnarray}
	where $z^{'}:=(z_{1},\cdots, z_{n})\in \mathbb{R}^{n}, y^{'}:=(y_{1},\cdots, y_{n})\in \mathbb{R}^{n}$.  Define
	\begin{eqnarray}
		\phi(y)=r_{0}y:=x,\notag
	\end{eqnarray}
	which maps $\bar{\mathbb{B}}^{n+1} $ to $\bar B_{r_{0}}$. Then map $\bar f:=\phi\circ f:\bar{\mathbb{R}}^{n+1}_+\to \bar{B}_{r_{0}}$  satisfies
	\begin{eqnarray}
		\bar{f}(\mathbb{R}^{n+1}_{+})=B_{r_{0}}, \quad \bar{f}(\partial \mathbb{R}^{n+1}_{+})=\partial B_{r_{0}}.\notag
	\end{eqnarray}
     $\bar{f}$ is a conformal transformation from $(\bar{\mathbb{R}}^{n+1}_{+}, \delta_{\mathbb{R}^{n+1}_{+}}:=|dz|^{2})$ to $(\bar{B}_{r_{0}}, \bar{g}_{\HH}=\frac{4}{(1-|x|^{2})^{2}}|dx|^{2})$ or $(\bar{B}_{r_{0}}, \bar{g}_{\SS}=\frac{4}{(1+|x|^{2})^{2}}|dx|^{2})$ (cf. \cite[Proposition 2.1]{QWX} ) and  
      	\begin{eqnarray}
      		\bar{f}^{\ast}\bar{g}_{\HH}=e^{2U}|dz|^{2},\quad \quad \bar{f}^{\ast}\bar{g}_{\SS}=e^{2V}|dz|^{2},\notag
      	\end{eqnarray}
      	where
      	\begin{eqnarray}
      		e^{U}:=\frac{4r_{0}}{(1-r_{0}^{2})(1+|z|^{2})+2(1+r_{0}^{2})z_{n+1}},\notag
      	\end{eqnarray} 
      and 
      \begin{eqnarray}
      		e^{V}:=\frac{4r_{0}}{(1+r_{0}^{2})(1+|z|^{2})+2(1-r_{0}^{2})z_{n+1}}.\notag
      \end{eqnarray}

	Let $\Sigma\subset \bar{B}_{r_{0}}$ be a properly embedded compact hypersurface with capillary boundary, given by an embedding $x:\bar{\mathbb{S}}^{n}_{+}\rightarrow \bar{B}_{r_{0}}$. We associate $\Sigma$ with a corresponding hypersurface $\widehat{\Sigma}\subset \bar{\mathbb{R}}^{n+1}_{+}$, given by embedding 
	\begin{eqnarray}
		\widehat{x}=\bar{f}^{-1}\circ x:~\bar{\mathbb{S}}^{n}_{+}\rightarrow \bar{\mathbb{R}}^{n+1}_{+}.\notag
	\end{eqnarray} 
	Since $(\bar{B}_{r_{0}}, \bar{g})$ and $(\bar{\mathbb{R}}^{n+1}_{+}, e^{2U}|dz|^{2})$ or $(\bar{\RR}^{n+1}_{+}, e^{2V}|dz|^{2})$ are isometric,  then $x: \bar{\mathbb{S}}^{n}_{+}\rightarrow \bar{B}_{r_{0}}$ can be identified as the embedding $\bar{x}: \bar{\mathbb{S}}^{n}_{+}\rightarrow (\bar{\mathbb{R}}^{n+1}_{+},  e^{2U}|dz|^{2})$ or $(\bar{\RR}^{n+1}_{+}, e^{2V}|dz|^{2})$.
	
	It is easy to check that $\widehat{X}_{a}:=(\bar{f}^{-1})_{\ast}(X_{a})=\frac{2r_{0}}{1+Kr_{0}^{2}}(\rho\partial_{\rho})$, 
	so a hypersurface $\widehat{\Sigma}\subset (\bar{\RR}^{n+1}_{+}, e^{2U}|dz|^{2})$ or  ($ \bar{\RR}^{n+1}_{+}, e^{2V}|dz|^{2} $)  is star-shaped (classically) with respect to the origin if and only if $\Sigma\subset (\bar{B}_{r_{0}},\bar{g})\subset M^{n+1}(-1)$ or $\Sigma\subset (\bar{B}_{r_{0}},\bar{g})\subset M^{n+1}(+1)$ is star-shaped with respect to $a\in \SS^n$. Then there exists some positive function $\rho(y)$ defined on $\bar{\mathbb{S}}_{+}^{n}$, such that
	\begin{eqnarray}
		\widehat{x}=\rho(y)y=\rho(\beta, \theta)y, \quad y:=(\beta,\theta)\in [0,\frac{\pi}{2}]\times \mathbb{S}^{n-1}.\notag
	\end{eqnarray}
	We denote $\text{div},\nabla$ as the divergence operator, covariant derivative on $\bar{\mathbb{S}}^{n}_{+}$ with respect to the standard spherical metric $\sigma$ on $\bar{\mathbb{S}}^{n}_{+}$ respectively. Set $u:=\log\rho$ and $v:=\sqrt{1+|\nabla u|^{2}}$. We have the following identities for $\HH^{n+1}$ and $\SS^{n+1}$ respectively.
	\begin{proposition}\label{basic formula hyperbolic}
		In a geodesic ball $B_{R}^{\HH}$ of $\HH^{n+1}$, it holds,
		\begin{enumerate}
			\item
			$ \bar{g}(X_{a}, \nu)=e^{U}\frac{2r_{0}}{1-r_{0}^{2}}\frac{\rho}{v}$.
			\item
			$V_{a}=\frac{2r_{0}(1-\rho^{2})}{(1-r_{0}^{2})(\rho^{2}+1)+2(1+r_{0}^{2})\rho\cos\beta}.$
			
			\item
			{\rm The mean curvature} $\widetilde{H}$ {\rm of} $\widehat{\Sigma}$ {\rm in} $(\bar{\mathbb{R}}^{n+1}_{+}, e^{2U}|dz|^{2})$ {\rm is}
			\begin{eqnarray*}
				\widetilde{H}&=&-\Bigg[\frac{1}{\rho v e^{U}}\sum\limits_{i,j=1}^{n}(\sigma^{ij}-\frac{u^{i}u^{j}}{v^{2}})u_{ij}+\frac{1+r_{0}^{2}}{2r_{0}}\frac{n\sin\beta \nabla_{\partial_{\beta}}u}{v}\\
				&&+\frac{n(\rho^{2}-1)(1-r_{0}^{2})}{4r_{0}\rho v}\Bigg].
			\end{eqnarray*}
			\item
			$\bar{g}(Y_{a},\nu)=-\frac{r_{0}\rho e^{U}}{v}+\frac{r_{0}^{2}-1}{2r_{0}}\frac{e^{U}}{2v}\left(\rho^{2}\cos\beta+2\rho+\cos\beta-(\rho^{2}-1)\sin\beta\nabla_{\partial_\beta}u\right).$
		\end{enumerate}
	\end{proposition}
	
	\begin{proof} (1)-(3) were  shown in \cite[Proposition 2.2]{QWX}. We only need to show (4). 
		Note that \begin{eqnarray*}   		\widehat{e}_{n+1}:&=&(\bar{f}^{-1})_{\ast}(e_{n+1})\\   		&=&\sum\limits_{i=1}^{n}\frac{\partial z_{i}}{\partial x_{n+1}}\frac{\partial}{\partial z_{i}}+\frac{\partial z_{n+1}}{\partial x_{n+1}}\frac{\partial}{\partial z_{n+1}}\\   		&=&\frac{1+z_{n+1}}{r_{0}}\sum\limits_{i=1}^{n}z_{i}\frac{\partial}{\partial z_{i}}+\frac{(1+z_{n+1})^{2}-|z|^{2}}{2r_{0}}\frac{\partial }{\partial z_{n+1}}			\\   		&=&\frac{1}{r_{0}}\rho\cos \beta(\rho \sin^{2}\beta\partial_{\rho}+\frac{\sin \beta}{2}\partial_{\beta})+\frac{1+\rho^{2}\cos2\beta}{2r_{0}}(\cos\beta\partial_{\rho}-\frac{\sin\beta}{\rho}\partial_{\beta})+\frac{1}{r_{0}}\rho\partial_{\rho}\\   		&=&\frac{\rho^{2}\cos\beta+2\rho+\cos\beta}{2r_{0}}\partial_{\rho}+\frac{(\rho^{2}-1)\sin\beta}{2r_{0}\rho}\partial_{\beta},  \end{eqnarray*} and
		the unit normal of $\widehat{\Sigma}\subset (\bar{\RR}_{+}^{n+1}, e^{2U}|dz|^{2})$ is    	\begin{eqnarray}   		\widehat{\nu}:=(\bar{f}^{-1})_{\ast}(\nu)=e^{-U}\frac{\partial_{\rho}-\rho^{-1}\nabla u}{v},\notag   	\end{eqnarray}
		By the definition $X_{a}$ and $Y_{a}$, we have   	\begin{eqnarray}   		Y_{a}=\frac{r_{0}^{2}-1}{2}e_{n+1}-\frac{1-r_{0}^{2}}{2}X_{a}.\notag   	\end{eqnarray} 
		Hence, it follows 
		\begin{eqnarray*}
			\bar{g}(Y_{a},\nu)&=& \bar g(\frac{r_{0}^{2}-1}{2}e_{n+1}-\frac{1-r_{0}^{2}}{2}X_{a},\nu)
			\\ 		&=&\bar{f}^{\ast}\bar{g}(\frac{r_{0}^{2}-1}{2}\widehat{e}_{n+1}-\frac{1-r_{0}^{2}}{2}\widehat{X}_{a},\widehat{\nu})\\
			&=&-\frac{r_{0}\rho e^{U}}{v}+\frac{(r_{0}^{2}-1)e^{U}}{2r_{0}}
			\left(\frac{\rho^{2}\cos\beta+2\rho+\cos\beta}{2v}-\frac{(\rho^{2}-1)\sin\beta}{2v}\nabla_{\partial_{\beta}}u\right).
		\end{eqnarray*}
	\end{proof}
Similarly, in the case $\SS^{n+1}$.
    \begin{proposition}\label{basic formula sphere}\
    	In a geodesic ball $B_{R}^{\SS}$ of $\SS^{n+1}$, it holds
    	\begin{enumerate}
    		\item
    		$ \bar{g}(X_{a}, \nu)=e^{V}\frac{2r_{0}}{1+r_{0}^{2}}\frac{\rho}{v}$.
    		\item
    		$V_{a}=\frac{2r_{0}(1-\rho^{2})}{(1+r_{0}^{2})(\rho^{2}+1)+2(1-r_{0}^{2})\rho\cos\beta}.$
    		
    		\item
    		{\rm The mean curvature} $\widetilde{H}$ {\rm of} $\widehat{\Sigma}$ {\rm in} $(\bar{\mathbb{R}}^{n+1}_{+}, e^{2V}|dz|^{2})$ {\rm is}
    		\begin{eqnarray*}
    			\widetilde{H}&=&-\Bigg[\frac{1}{\rho v e^{V}}\sum\limits_{i,j=1}^{n}(\sigma^{ij}-\frac{u^{i}u^{j}}{v^{2}})u_{ij}+\frac{1-r_{0}^{2}}{2r_{0}}\frac{n\sin\beta \nabla_{\partial_{\beta}}u}{v}\\
    			&&+\frac{n(\rho^{2}-1)(1+r_{0}^{2})}{4r_{0}\rho v}\Bigg].
    		\end{eqnarray*}
    		\item
    		$\bar{g}(Y_{a},\nu)=\frac{r_{0}\rho e^{V}}{v}-\frac{r_{0}^{2}+1}{2r_{0}}\frac{e^{V}}{2v}\left(\rho^{2}\cos\beta+2\rho+\cos\beta-(\rho^{2}-1)\sin\beta\nabla_{\partial_\beta}u\right).$
    	\end{enumerate}
    \end{proposition}
The proof of Proposition \ref{basic formula sphere} is the similar to Proposition \ref{basic formula hyperbolic}, we omit it here.

With the help of Proposition  \ref{basic formula hyperbolic}  and \ref{basic formula sphere},  following the argument as in  \cite [Section 3.3]{WW20}, we can reduce the first equation in \eqref{flow with capillary in hyperbolic} to the following scalar equation
	\begin{eqnarray}\label{saclar equation in hyperbolic}
		\partial_{t}u=\frac{v}{\rho e^{U}}\widehat{F} .
	\end{eqnarray} 	
Moreover, in the case $\HH^{n+1}$, 
	\begin{eqnarray*}
		\widehat{F} &:=&nV_{a}+n\sinh R\cos\theta\bar{g}(Y_{a},\nu)-\widetilde{H}\bar{g}(X_{a},\nu)\\
		&=&\frac{2nr_{0}(1-\rho^{2})}{(1-r_{0}^{2})(1+\rho^{2})+2(1+r_{0}^{2})\rho\cos\beta}\frac{|\nabla u|^{2}}{v^{2}}-\frac{2n\cos\theta r_{0}^{2}}{1-r_{0}^{2}}\frac{\rho e^{U}}{v}\\
		&&-\frac{n\cos\theta e^{U}}{2v}\left(\rho^{2}\cos\beta+2\rho+\cos\beta-(\rho^{2}-1)\sin\beta\nabla_{\partial_{\beta}}u\right)\\
		&&+
		\frac{n(1+r_{0}^{2})\rho e^{U}\sin\beta\nabla_{\partial_{\beta}}u}{(1-r_{0}^{2})v^{2}}+\frac{2r_{0}}{1-r_{0}^{2}}\frac{1}{v^{2}}\sum\limits_{i,j=1}^{n}(\sigma^{ij}-\frac{u^{i}u^{j}}{v^{2}})u_{ij},
	\end{eqnarray*}
and it is now easy to see that \eqref{saclar equation in hyperbolic} is equivalent to
	\begin{eqnarray}\notag
		\partial_{t}u
		&=&\frac{2r_{0}}{1-r_{0}^{2}}{\rm div}(\frac{\nabla u}{\rho v e^{U}})-\frac{2(n+1)r_{0}}{(1-r_{0}^{2})v}\sigma \left(\nabla u,\nabla(\frac{1}{\rho e^{U}})\right)-\frac{2n\cos\theta r_{0}^{2}}{1-r_{0}^{2}}\\
	\label{G hyperbolic}	&&-\frac{n\cos\theta}{2\rho}\left(\rho^{2}\cos\beta+2\rho+\cos\beta-(\rho^{2}-1)\sin\beta\nabla_{\partial_{\beta}}u\right)\\ 
		&:=&G(\nabla^{2}u,\nabla u,\rho, \beta).\notag	
	\end{eqnarray}

   Similarly, in the case $\SS^{n+1}$,  $\widehat{F}$ in flow \eqref{saclar equation in hyperbolic} has the form 
\begin{eqnarray*}
	\widehat{F} &:=&nV_{a}+n\sin R\cos\theta\bar{g}(Y_{a},\nu)-\widetilde{H}\bar{g}(X_{a},\nu)\\
	&=&\frac{2nr_{0}(1-\rho^{2})}{(1+r_{0}^{2})(1+\rho^{2})+2(1-r_{0}^{2})\rho\cos\beta}\frac{|\nabla u|^{2}}{v^{2}}+\frac{2n\cos\theta r_{0}^{2}}{1+r_{0}^{2}}\frac{\rho e^{V}}{v}\\
	&&-\frac{n\cos\theta e^{V}}{2v}\left(\rho^{2}\cos\beta+2\rho+\cos\beta-(\rho^{2}-1)\sin\beta\nabla_{\partial_{\beta}}u\right)\\
	&&+
	\frac{n(1-r_{0}^{2})\rho e^{V}\sin\beta\nabla_{\partial_{\beta}}u}{(1+r_{0}^{2})v^{2}}+\frac{2r_{0}}{1+r_{0}^{2}}\frac{1}{v^{2}}\sum\limits_{i,j=1}^{n}(\sigma^{ij}-\frac{u^{i}u^{j}}{v^{2}})u_{ij},
\end{eqnarray*}which follows that \eqref{saclar equation in hyperbolic} is equivalent to 
\begin{eqnarray}	\notag
	\partial_{t}u
     &=&\frac{2r_{0}}{1+r_{0}^{2}}{\rm div}(\frac{\nabla u}{\rho v e^{V}})-\frac{2(n+1)r_{0}}{(1+r_{0}^{2})v}\sigma \left(\nabla u,\nabla(\frac{1}{\rho e^{V}})\right)+\frac{2n\cos\theta r_{0}^{2}}{1+r_{0}^{2}}\\\label{G sphere}
	&&-\frac{n\cos\theta}{2\rho}\left(\rho^{2}\cos\beta+2\rho+\cos\beta-(\rho^{2}-1)\sin\beta\nabla_{\partial_{\beta}}u\right)\\
	&:=&G(\nabla^{2}u,\nabla u,\rho, \beta). 	\notag
\end{eqnarray}
	
Next we derive the boundary condition. We show the case $\HH^{n+1}$ here, since the case $\SS^{n+1}$ is similar.  The capillary boundary condition in flow \eqref{flow with capillary in hyperbolic} implies 
	\begin{eqnarray*}
		-\cos\theta =\bar{g}(\nu,\bar{N}\circ x)=e^{2U}\left\<\widehat{\nu}, (\bar{f}^{-1})_{\ast}(\bar{N}\circ x)\right\>.
	\end{eqnarray*}
	Note that $(\bar{f}^{-1})_{\ast}(\bar{N}\circ x)=-e^{-U}\frac{\partial_{\beta}}{\rho}$ on $\partial \mathbb{R}^{n+1}_{+}$, then
	\begin{eqnarray}\label{boundary condition}
		\nabla_{\partial_{\beta}}u=\cos\theta \sqrt{1+|\nabla u|^{2}}, \quad {\rm on}~\partial\mathbb{S}^{n}_{+}.
	\end{eqnarray}
  
  In summary, the flow \eqref{flow with capillary in hyperbolic} is equivalent to the following scalar parabolic equation on $\bar{\mathbb{S}}^{n}_{+}$,	 
	\begin{eqnarray}\label{scalar equ in hyperbolic}
		\begin{cases}
			\frac{\partial  u}{\partial t} =G (\nabla^2 u,\nabla u, \rho ,\beta)\quad &\text{in}  \quad \mathbb{S}^n_+\times[0,T),\\
			\nabla_{\partial_{\beta}}u=\cos\theta \sqrt{1+|\n u|^2}  \quad & \text{on}  \quad \partial\mathbb{S}^n_+\times[0,T),\\
			u(\cdot,0)=u_0(\cdot) \quad &\text{in}  \quad \mathbb{S}^n_+,
		\end{cases}
	\end{eqnarray}
   where $G$ has the form as \eqref{G hyperbolic} in $\HH^{n+1}$ and \eqref{G sphere} in $\SS^{n+1}$ respectively, $u_{0}:=\log \rho_{0}$ and $\rho_{0}$ corresponds to the hypersurface $x_{0}(M)$ under the transformation $\bar f$.

	\section{A priori estimates and convergence}\label{sec3}
	In this section, we focus on establishing the uniform height and gradient estimates for the solution of scalar parabolic equation \eqref{scalar equ in hyperbolic}. 	The key ingredient is gradient estimate. Since the expression of $G$ in \eqref{scalar equ in hyperbolic} for $\HH^{n+1}$ and $\SS^{n+1}$  are  similar, we only show the case $\HH^{n+1}$ below. One can obtain the same a priori estimates for the case $\SS^{n+1}$ by just adapting the same approach as below with minor modifications.

	The short-time existence of the flow \eqref{flow with capillary in hyperbolic}  follows from the standard  PDE theory (cf. \cite{HP1999}), due to our assumption  of star-shaped, $\bar{g}(X_{a}, \nu)>0$ on $x_0(M)$.  Next, in order to establish the long-time existence of the flows, we need to obtain the uniform height and  gradient estimates for the solutions of flow , then the long-time existence and  uniform $C^{\infty}$ estimates follows from the standard quasi-linear parabolic PDE theory with strictly oblique boundary condition (cf. \cite{LSU}).

	In the following, we use the Einstein summation convention, i.e., if not stated otherwise, the repeated arabic indices $i,j,k$ should be summed from $1$ to $n$.
	Besides, we introduce the following notations. $a^{ij}:=(\sigma^{ij}-\frac{u^{i}u^{j}}{v^{2}})$,   $u_{\beta}:=\nabla_{\partial_{\beta}}u=\s(\n u,\p_\beta)$ and
	\begin{eqnarray*}
		G^{ij}&:=&\frac{\partial G(r, p,\rho, \beta)}{\partial r_{ij}}\Big|_{r=\nabla^{2}u, p=\nabla u}=\frac{2r_{0}}{1-r_{0}^{2}}\frac{1}{\rho v e^{U}}a^{ij},\\
		G_{p_{k}}&:=&\frac{\partial G(r, p, \rho,\beta)}{\partial p_{k}}\Big|_{r=\nabla^{2}u, p=\nabla u}\\
		&=&-\frac{2r_{0}}{1-r_{0}^{2}}\frac{u_{k}}{e^{U}\rho v^{3}}a^{ij}u_{ij}
		+\frac{2r_{0}}{1-r_{0}^{2}}\frac{1}{e^{U}\rho v}(-\frac{\delta_{i}^{k}u^{j}+\delta_{j}^{k}u^{i}}{v^{2}}
		+\frac{2u_{k}u^{i}u^{j}}{v^{4}})u_{ij}\\
		&&+\frac{n(1+r_{0}^{2})}{1-r_{0}^{2}}\frac{\sin\beta\sigma(\partial_{\beta},\partial_{k})}{v}-\frac{n(1+r_{0}^{2})}{1-r_{0}^{2}}\frac{\sin\beta u_{\beta} u_{k}}{v^{3}}
		+\frac{n(1-\rho^{2})}{\rho}\frac{u_{k}}{v}\\
		&&+\frac{n(\rho^{2}-1)|\nabla u|^{2}u_{k}}{2\rho v^{3}}+\frac{\rho^{2}-1}{2\rho}\sin\beta \sigma(\partial_{\beta},\partial_{k}),\\
	    G_{\rho}&:=&\frac{\partial G (r, p,\rho, \beta)}{\partial \rho}\Big|_{r=\nabla^{2}u, p=\nabla u}\\
		&=&\frac{1}{2v}(1-\frac{1}{\rho^{2}})a^{ij}u_{ij}-\frac{n}{2}(1+\frac{1}{\rho^{2}})\frac{|\nabla u|^{2}}{v}
		-\frac{n\cos\theta\cos\beta}{2}(1-\frac{1}{\rho^{2}})\\
		&&+\frac{n}{2}(1+\frac{1}{\rho^{2}})\cos\theta\sin\beta u_{\beta},
	\end{eqnarray*}
   \begin{eqnarray*}
		G_{\beta}&:=&\frac{\partial G(r, p,\rho, \beta)}{\partial \beta}\Big|_{r=\nabla^{2}u, p=\nabla u}\\
		&=&\frac{1+r_{0}^{2}}{1-r_{0}^{2}}\frac{-a^{ij}u_{ij}\sin\beta+n\cos\beta u_{\beta}}{v}+\frac{n\cos\theta}{2\rho}\Big[\sin\beta \rho^{2}+\sin\beta\\
		&&+(\rho^{2}-1)\cos\beta u_{\beta}\Big].	
	\end{eqnarray*}
	Firstly, we have the height estimate for the solution of flow \eqref{scalar equ in hyperbolic} for the case $\HH^{n+1}$.
	\begin{proposition}\label{pro3.2}
		Assume   $\Sigma_0:=x_{0}(M)\subset B_R^{\HH}$ is star-shaped with respect to $a\in \SS^n$ and satisfies
		\begin{eqnarray}\label{star shaped}
			\Sigma_0\subset C_{\theta, r_{1}}(a) \setminus C_{\theta, r_{2}}(a), 
		\end{eqnarray}
		for some $0<r_{1}<r_{2}$, where $C_{\theta, r}(a)$ is defined by \eqref{static model}. Then the solution $\Sigma_{t}:=x(M,t)$ of \eqref{flow with capillary in hyperbolic} satisfies
		\begin{eqnarray}
			\Sigma_{t}\subset  C_{\theta, r_{2}}(a)\setminus C_{\theta, r_{1}}(a).\notag
		\end{eqnarray}
		Moreover, if $u $ solves \eqref{scalar equ in hyperbolic} and $G$ has the form \eqref{G hyperbolic}, then
		\begin{eqnarray}
			||u||_{C^{0}\left(\bar{\mathbb{S}}^{n}_{+}\times [0,T)\right)}\leq C,\notag
		\end{eqnarray} 
		where $C$ is a positive constant, depending only on the initial datum.
	\end{proposition}
	\begin{proof}
		Since the spherical cap $C_{\theta, r}(a)$ is the static solution to flow \eqref{flow with capillary in hyperbolic} for each $r>0$, that is, it satisfies
		\begin{eqnarray}
			nV_{a}+n\sinh R\cos\theta \bar{g}(Y_{a}, \nu)-H\bar{g}(X_{a}, \nu)=0.\notag
		\end{eqnarray}
		Then the assertion follows from the avoidance principle  for the strictly parabolic equation with a capillary boundary condition (see 
		\cite [Proposition 4.2]{WW20}).
	\end{proof}

	In order to get the gradient estimate, we need to employ the distance function $d(x):=dist_{\sigma}(x,\partial\mathbb{S}^{n}_{+})$. It is well-known that $d$ is well-defined and smooth for $x$ near $\partial\mathbb{S}^{n}_{+}$ and $\nabla d=-\partial_{\beta}$ on $\partial \mathbb{S}^{n}_{+}$, where $\partial_{\beta}$ is the unit outer normal vector field of $\partial \mathbb{S}^{n}_{+}$ in $\SS^n_+$. We can extend $d$ to be a smooth function defined in $\bar{\mathbb{S}}^{n}_{+}$ and satisfy that
	$$d\geq 0,\quad\quad |\nabla d|\leq1 ,\quad  \text{  in } \bar{\mathbb{S}}^{n}_{+}.$$

	We will use $O(s)$ to denote terms that are bounded by  $Cs$ for some constant $C>0$, which    depends only on the $C^{0}$ norm of $u$. And the constant $C$ may change from line to line. Next, we choose a suitable auxiliary function that had been used in \cite{WW20} to  obtain the uniform gradient estimate for the flow \eqref{scalar equ in hyperbolic} in the case $\HH^{n+1}$. 
	\begin{proposition}\label{pro3.3}
		If $u:\bar{\mathbb{S}}^{n}_{+}\times[0,T)\rightarrow \mathbb{R}$ solves  \eqref{scalar equ in hyperbolic} and $G$ has the form \eqref{G hyperbolic},  $|\cos\theta|<\frac{3n+1}{5n-1}$, then for any $(x,t)\in \bar{\mathbb{S}}^{n}_{+}\times [0,T)$, 
		\begin{eqnarray}
			|\nabla u|(x,t)\leq  C,\notag
		\end{eqnarray}
		where $C$ is a positive constant, depending only on the initial datum.
	\end{proposition}	
	\begin{proof}
		Define the function as
		\begin{eqnarray}
			\Phi:=(1+Kd)v+\cos\theta \sigma(\nabla u,\nabla d),\notag
		\end{eqnarray}
		where $K$ is a positive constant to be determined later. For any $T'<T$, assume $\Phi$ attains its maximum value at some point, say $(x_{0}, t_{0})\in\bar{\mathbb{S}}^{n}_{+}\times [0,T']$. 
		
		Following the same argument as in \cite[Proposition 4.3, Case 1]{WW20}, by choosing  $K>0$ sufficiently large,  $\Phi$ does not attian its maximum value on $\partial \mathbb{S}^{n}_{+}$, hence we have either $x_0\in \SS^n_+$ or $t_0=0$. If   $t_0=0$, it is easy to see
		\begin{eqnarray}
			\sup\limits_{\bar{\mathbb{S}}^{n}_{+}\times [0,T']}  |\n u|\leq C,
		\end{eqnarray} 
		where $C$ is a positive constant depending only on $n$ and $u_{0}$. 
		
		Next we analyze the case $x_0\in\SS^n_+$ and complete the gradient estimate.   By rotating the geodesic coordinate $\{\frac{\partial}{\partial x_{i}}\}_{i=1}^{n}$ at $x_0$, we assume that
		$$|\nabla u|=u_{1}>0,~{\rm and} ~\{u_{\alpha\beta}\}_{2\leq \alpha,\beta\leq n} ~{\rm is~diagonal}.$$
		Assume that $u_{1}(x_{0}, t_{0})$ is large enough, otherwise we finish the proof. All the computation below are done at the point   $(x_{0}, t_{0})$. Note that
		\begin{eqnarray*}
			0=\nabla_{i}\Phi= (1+Kd)v_{i}+Kd_{i}v+\cos\theta(u_{li}d_{l}+u_{l}d_{li}), ~1\leq i\leq n.
		\end{eqnarray*}
		It follows that
		\begin{eqnarray}
			\left[(1+Kd)\frac{u_{1}}{v}+\cos\theta d_{1}\right]u_{1\alpha}=\cos\theta u_{\alpha\alpha}d_{\alpha}-\cos\theta u_{a}d_{1\alpha}-Kd_{\alpha} v,\notag
		\end{eqnarray}
		and 
		\begin{eqnarray}
			\left[(1+Kd)\frac{u_{1}}{v}+\cos\theta d_{1}\right]u_{11}=-\cos\theta u_{\alpha1} d_{\alpha}-\cos\theta u_{1}d_{11}-Kd_{1}v.\notag
		\end{eqnarray}
		Denote $S:=(1+Kd)\frac{u_{1}}{v}+\cos\theta d_{1}$, it is easy to see that $0<C(\delta, \theta)\leq S\leq 2+K$, if we assume $u_{1}\geq \delta>0$, otherwise we complete the proof. Hence
		\begin{eqnarray}\label{3.2}
			u_{1\alpha}&=&-\frac{\cos\theta d_{\alpha}}{S}u_{\alpha\alpha}-\frac{1}{S}(\cos\theta u_{1}d_{1\alpha}+Kd_{\alpha}v)\notag\\
			&=&-\frac{\cos\theta d_{\alpha}}{S}u_{\alpha\alpha}+O(v), ~2\leq \alpha\leq n,
		\end{eqnarray}
		and
		\begin{eqnarray}\label{3.3}
			u_{11}&=&-\frac{1}{S}\cos\theta u_{\alpha 1}d_{\alpha}+\frac{1}{S}(-\cos\theta u_{1}d_{11}-Kd_{1}v)\notag\\
			&=&\frac{\cos^{2}\theta}{S^{2}}\sum\limits_{\alpha=2}^{n}d_{\alpha}^{2}u_{\alpha\alpha}+O(v).
		\end{eqnarray}
		On the other hand, we have
		\begin{eqnarray} \notag
			0&\leq &(\partial_{t}-G^{ij}\nabla_{ij}-G_{p_{i}}\nabla_{i})\Phi \\ \notag
			&=&\frac{(1+Kd)}{v}u_{l}(u_{lt}-G^{ij}u_{lij}-G_{p_{i}}u_{li})+d_{k}\cos\theta (u_{kt}-G^{ij}u_{kij}-G_{p_{i}}u_{ki}) \\ \notag
			&&+(1+Kd)(\frac{G^{ij}u_{l}u_{li}u_{k}u_{kj}}{v^{3}}-\frac{G^{ij}u_{li}u_{lj}}{v})-(2G^{ij}u_{ki}d_{kj}\cos\theta+2KG^{ij}d_{i}v_{j})\\ \notag
			&&-(G^{ij}u_{k}d_{kij}\cos\theta+KG^{ij}d_{ij}v)-G_{p_{i}}(Kd_{i}v+\cos\theta u_{k}d_{ki})\\  
			& 	:=&I_{1}+I_{2}+I_{3}+I_{4}+I_{5}+I_{6}.\label{3.4.}
		\end{eqnarray}
		Next, we carefully handle those terms one by one. Differentiating the first equation in \eqref{scalar equ in hyperbolic}, 
		\begin{eqnarray}
			u_{tl}=G^{ij}u_{ijl}+G_{p_{i}}u_{il}+G_{\rho}\rho u_{l}+G_{\beta}\sigma(\partial_{\beta},\partial_{l}),\notag
		\end{eqnarray}
		Combining with the communicative formula on $\mathbb{S}^{n}_{+}$,
		\begin{eqnarray}
			u_{ijl}=u_{lij}+u_{j}\s_{li}-u_{l}\sigma_{ij},\notag
		\end{eqnarray}
		it follows
		\begin{eqnarray*}
			u_{tl}=G^{ij}u_{lij}+G^{ij}u_{j}\sigma_{li}-\sum\limits_{i=1}^{n}G^{ii}u_{l}+G_{p_{i}}u_{il}+G_{\rho}\rho u_{l}+G_{\beta}\sigma(\partial_{\beta},\partial_{l}).	
		\end{eqnarray*}
		First, we deal with the term $I_{1}$. 
		\begin{eqnarray*}
			I_{1}&=&\frac{(1+Kd)}{v}u_{l}(u_{lt}-G^{ij}u_{lij}-G_{p_{i}}u_{il})\\
			&=&\frac{(1+Kd)}{v}G^{ij}u_{j}u_{l}\sigma_{li}
			+\frac{(1+Kd)u_{l}}{v}\left(G_{\rho}\rho u_{l}+G_{\beta}\sigma(\partial_{\beta},\partial_{l})\right)\\
			&&-\frac{(1+Kd)|\nabla u|^{2}}{v}\sum\limits_{i=1}^{n}G^{ii}\\
			&=&\left[\frac{(1+Kd)|\nabla u|^{2}}{2v^{4}}(\rho-\frac{1}{\rho})u_{11}-
			\frac{1+r_{0}^{2}}{1-r_{0}^{2}}\frac{(1+Kd)\sin\beta u_{\beta}}{v^{4}}u_{11}\right]\\
			&&+\left[\frac{(1+Kd)|\nabla u|^{2}}{2v^{2}}(\rho-\frac{1}{\rho})\sum\limits_{\alpha=2}^{n}u_{\alpha\alpha}\right]\\
			&&-\left[\frac{n(1+Kd)|\nabla u|^{4}}{2v^{2}}(\rho+\frac{1}{\rho})-\frac{n(1+Kd)|\nabla u|^{2}}{2v}(\rho+\frac{1}{\rho})\cos\theta\sin\beta u_{\beta}\right]\\
			&&-\Bigg[\frac{1+r_{0}^{2}}{1-r_{0}^{2}}\frac{(1+Kd)u_{\beta}\sin\beta}{v^{2}}\sum\limits_{\alpha=2}^{n}u_{\alpha\alpha}+\frac{1+r_{0}^{2}}{1-r_{0}^{2}}\frac{n(1+Kd)\cos\beta u_{\beta}^{2}}{v^{2}}\\
			&&-\frac{n\cos\theta(1+Kd)u_{\beta}}{2\rho v}\left(\sin\beta\rho^{2}+\sin\beta +(\rho^{2}-1)\cos\beta u_{\beta}\right)\\
			&&-\frac{2r_{0}}{1-r_{0}^{2}}\frac{(1+Kd)(1-n)|\nabla u|^{2}}{\rho v^{2} e^{U}}
			+\frac{n(1+Kd)|\nabla u|^{2}\cos\beta\cos\theta}{2v}(\rho-\frac{1}{\rho})\Bigg]\\
			&:=&I_{11}+I_{12}+I_{13}+I_{14}.
		\end{eqnarray*}		
		By \eqref{3.3}, we obtain
		\begin{eqnarray*}
			I_{11}&=&\frac{(1+Kd)|\nabla u|^{2}}{2v^{4}}(\rho-\frac{1}{\rho})u_{11}-
			\frac{1+r_{0}^{2}}{1-r_{0}^{2}}\frac{(1+Kd)\sin\beta u_{\beta}}{v^{4}}u_{11}\\
			&=&\left[\frac{(1+Kd)|\nabla u|^{2}}{2v^{4}}(\rho-\frac{1}{\rho})-
			\frac{1+r_{0}^{2}}{1-r_{0}^{2}}\frac{(1+Kd)\sin\beta u_{\beta}}{v^{4}}\right]\\
			&&\cdot \left(\frac{\cos^{2}\theta}{S^{2}}\sum\limits_{\alpha=2}^{n}d_{\alpha}^{2}u_{\alpha\alpha}+O(v)\right)\\
			&=&O(\frac{1}{v^{2}})\sum\limits_{\alpha=2}^{n}|u_{\alpha\alpha}|+O(\frac{1}{v}).
		\end{eqnarray*}
		and
		\begin{eqnarray*}
			I_{14}&=&-\Bigg[\frac{1+r_{0}^{2}}{1-r_{0}^{2}}\frac{(1+Kd)u_{\beta}\sin\beta}{v^{2}}\sum\limits_{\alpha=2}^{n}u_{\alpha\alpha}+\frac{1+r_{0}^{2}}{1-r_{0}^{2}}\frac{n(1+Kd)\cos\beta u_{\beta}^{2}}{v^{2}}\\
			&&-\frac{n\cos\theta(1+Kd)u_{\beta}}{2\rho v}\left(\sin\beta\rho^{2}+\sin\beta +(\rho^{2}-1)\cos\beta u_{\beta}\right)\\
			&&-\frac{2r_{0}}{1-r_{0}^{2}}\frac{(1+Kd)(1-n)|\nabla u|^{2}}{\rho v^{2} e^{U}}
			+\frac{n(1+Kd)|\nabla u|^{2}\cos\beta\cos\theta}{2v}(\rho-\frac{1}{\rho})\Bigg]\\
			&&=O(\frac{1}{v^{2}})\sum\limits_{\alpha=2}^{n}|u_{\alpha\alpha}|+O(v).	
		\end{eqnarray*}
		For the term $I_{2}$,
		\begin{eqnarray*}
			I_{2}&=&d_{k}\cos\theta (u_{kt}-G^{ij}u_{kij}-G_{p_{i}}u_{ki})\\
			&=&\cos\theta \sigma(\nabla u, \nabla d)\rho G_{\rho}+\cos\theta G_{\beta} d_{\beta}+
			\cos\theta G^{11}\sigma(\nabla u,\nabla d)\\
			&&-\cos\theta \sigma(\nabla u,\nabla d)\sum\limits_{i=1}^{n}G^{ii}\\
			&=&\left[\frac{\cos\theta\sigma(\nabla u,\nabla d)}{2v^{3}}(\rho-\frac{1}{\rho})u_{11}-\frac{1+r_{0}^{2}}{1-r_{0}^{2}}\frac{\sin\beta\cos\theta d_{\beta}}{v^{3}}u_{11}\right]\\
			&&+\left[\frac{\cos\theta \sigma(\nabla u,\nabla d)}{2v}(\rho-\frac{1}{\rho})\sum\limits_{\alpha=2}^{n}u_{\alpha\alpha}\right]\\
			&&-\left[\frac{n\cos\theta \sigma(\nabla u,\nabla d)}{2}
			(\rho+\frac{1}{\rho})(\frac{|\nabla u|^{2}}{v}-\cos\theta\sin\beta \nabla_{\partial_{\beta}}u)\right]\\
			&&-\Bigg[\frac{n\cos^{2}\theta\cos\beta \sigma(\nabla u,\nabla d)}{2}(\rho-\frac{1}{\rho})-\frac{n\cos\theta\cos\beta d_{\beta}\nabla_{\partial_{\beta}}u}{v}\\
			&&+\frac{2(n-1)r_{0}}{\rho v e^{U}(1-r_{0}^{2})}\cos\theta \sigma(\nabla u,\nabla d)+\frac{1+r_{0}^{2}}{1-r_{0}^{2}}\frac{\cos\theta\sin\beta d_{\beta}}{v}\sum\limits_{\alpha=2}^{n}u_{\alpha\alpha}\\
			&&-\frac{n\cos^{2}\theta d_{\beta}}{2\rho}\left(\sin\beta\rho^{2}+\sin\beta+(\rho^{2}-1)\cos\beta\nabla_{\partial_{\beta}}u\right)\Bigg]\\
			&:=&I_{21}+I_{22}+I_{23}+I_{24}.
		\end{eqnarray*}
		For the term $I_{21}, I_{24}$, we see
		\begin{eqnarray*}
			I_{21}&=&\frac{\cos\theta\sigma(\nabla u,\nabla d)}{2v^{3}}(\rho-\frac{1}{\rho})u_{11}-\frac{1+r_{0}^{2}}{1-r_{0}^{2}}\frac{\sin\beta\cos\theta d_{\beta}}{v^{3}}u_{11}\\
			&=&\left[\frac{\cos\theta\sigma(\nabla u,\nabla d)}{2v^{3}}(\rho-\frac{1}{\rho})-\frac{1+r_{0}^{2}}{1-r_{0}^{2}}\frac{\sin\beta\cos\theta d_{\beta}}{v^{3}}\right]\\
			&&\cdot \left(\frac{\cos^{2}\theta}{S^{2}}\sum\limits_{\alpha=2}^{n}d_{\alpha}^{2}u_{\alpha\alpha}+O(v)\right)\\
			&=&O(\frac{1}{v^{2}})\sum\limits_{\alpha=2}^{n}|u_{\alpha\alpha}|+O(\frac{1}{v^{2}}),
		\end{eqnarray*}
		and 
		\begin{eqnarray*}
			I_{24}&=&-\Bigg[\frac{n\cos\theta^{2}\cos\beta \sigma(\nabla u,\nabla d)}{2}(\rho-\frac{1}{\rho})-\frac{n\cos\theta\cos\beta d_{\beta}\nabla_{\partial_{\beta}}u}{v}\\
			&&+\frac{2(n-1)r_{0}}{\rho v e^{U}(1-r_{0}^{2})}\cos\theta \sigma(\nabla u,\nabla d)+\frac{1+r_{0}^{2}}{1-r_{0}^{2}}\frac{\cos\theta\sin\beta d_{\beta}}{v}\sum\limits_{\alpha=2}^{n}u_{\alpha\alpha}\\
			&&-\frac{n\cos^{2}\theta d_{\beta}}{2\rho}\left(\sin\beta\rho^{2}+\sin\beta+(\rho^{2}-1)\cos\beta\nabla_{\partial_{\beta}}u\right)\Bigg]\\
			&=&O(\frac{1}{v})\sum\limits_{\alpha=2}^{n}|u_{\alpha\alpha}|+O(v).
		\end{eqnarray*}
		Next, we handle the term $I_{3}$.
		\begin{eqnarray*}
			I_{3}&=&(1+Kd)(\frac{G^{ij}u_{l}u_{li}u_{k}u_{ki}}{v^{3}}-\frac{G^{ij}u_{li}u_{lj}}{v})\\
			&=&\frac{2r_{0}}{1-r_{0}^{2}}\frac{(1+Kd)}{\rho v e^{U}}(-\frac{1}{v^{5}}u_{11}^{2}-\frac{2}{v^{3}}\sum\limits_{\alpha=2}^{n}u_{1\alpha}^{2})
			-\frac{2r_{0}}{1-r_{0}^{2}}\frac{(1+Kd)}{\rho v^{2}e^{U}}\sum\limits_{\alpha=2}^{n}u_{\alpha\alpha}^{2}\\
			&=&\frac{2r_{0}}{1-r_{0}^{2}}\frac{(1+Kd)}{\rho v e^{U}}(-\frac{1}{v^{5}}u_{11}^{2}-\frac{2}{v^{3}}\sum\limits_{\alpha=2}^{n}u_{1\alpha}^{2})\\
			&&-(1-\varepsilon)\frac{2r_{0}}{1-r_{0}^{2}}\frac{(1+Kd)}{\rho v^{2}e^{U}}\sum\limits_{\alpha=2}^{n}u_{\alpha\alpha}^{2}-\varepsilon\frac{2r_{0}}{1-r_{0}^{2}}\frac{(1+Kd)}{\rho v^{2}e^{U}}\sum\limits_{\alpha=2}^{n}u_{\alpha\alpha}^{2}\\
			&:=&I_{31}+I_{32}+I_{33}.
		\end{eqnarray*}
		Hence
		\begin{eqnarray*}
			I_{31}&=&\frac{2r_{0}}{1-r_{0}^{2}}\frac{(1+Kd)}{\rho v e^{U}}(-\frac{1}{v^{5}}u_{11}^{2}-\frac{2}{v^{3}}\sum\limits_{\alpha=2}^{n}u_{1\alpha}^{2})\\
			&\leq& O(\frac{1}{v^{4}})\sum\limits_{\alpha=2}^{n}u_{\alpha\alpha}^{2}+O(v).
		\end{eqnarray*}
		Finally, we deal with the other remaining  terms in \eqref{3.4.}.
		\begin{eqnarray*}
			I_{4}+I_{5}+I_{6}&=&-(2G^{ij}u_{ki}d_{kj}\cos\theta+2KG^{ij}d_{i}v_{j})
			-(G^{ij}u_{k}d_{kij}\cos\theta+KG^{ij}d_{ij}v)\\
			&&-G_{p_{i}}(Kd_{i}v+\cos\theta u_{k}d_{ki})\\
			&=&O(\frac{1}{v})\sum\limits_{\alpha=2}^{n}|u_{\alpha\alpha}|+O(v).
		\end{eqnarray*}
		By the arithmetic-geometric inequality, we have
		\begin{eqnarray*}
			I_{12}+I_{22}+I_{32}&=&\frac{(1+Kd)|\nabla u|^{2}}{2v^{2}}(\rho-\frac{1}{\rho})\sum\limits_{\alpha=2}^{n}u_{\alpha\alpha}+\frac{\cos\theta \sigma(\nabla u,\nabla d)}{2v}(\rho-\frac{1}{\rho})\sum\limits_{\alpha=2}^{n}u_{\alpha\alpha}\\
			&&-(1-\varepsilon)\frac{2r_{0}}{1-r_{0}^{2}}\frac{(1+Kd)}{\rho v^{2}e^{U}}\sum\limits_{\alpha=2}^n u_{\alpha\alpha}^{2}\\
			&=&\frac{u_{1}}{2v}S(\rho-\frac{1}{\rho})\sum\limits_{\alpha=2}^{n}u_{\alpha\alpha}-(1-\varepsilon)\frac{2r_{0}}{1-r_{0}^{2}}\frac{(1+Kd)}{\rho v^{2}e^{U}}\sum\limits_{\alpha=2}^n u_{\alpha\alpha}^{2}\\
			&\leq&\frac{1-r_{0}^{2}}{2r_{0}}\frac{(n-1)}{1+Kd}\frac{(\rho-\frac{1}{\rho})^{2}S^{2}\rho e^{U}}{16(1-\varepsilon)}u_{1}^{2}\\
			&\leq &\frac{1-r_{0}^{2}}{2r_{0}}\frac{(n-1)(1+|\cos\theta|)S}{16(1-\varepsilon)}(\rho-\frac{1}{\rho})^{2}\rho e^{U}u_{1}^{2}.
		\end{eqnarray*}
		Fix a positive constant $a_{0}\in \left(|\cos\theta|, \frac{3n+1}{5n-1}\right)$, if
		\begin{eqnarray}
			\frac{|\nabla u|^{2}}{v}-\cos\theta\sin\beta u_{\beta}< (1-a_{0})u_{1},\notag
		\end{eqnarray}
		then
		\begin{eqnarray}
			\begin{aligned}
				\frac{|\nabla u|^{2}}{v}-|\cos\theta|u_{1}&\leq \frac{|\nabla u|^{2}}{v}-\cos\theta\sin\beta u_{\beta}<(1-a_{0})u_{1},\notag
			\end{aligned}
		\end{eqnarray}
		which implies
		\begin{eqnarray*}
			u_{1}^{2}\leq \frac{\left(\cos\theta+(1-a_{0})\right)^{2}}{1-\left(\cos\theta+(1-a_{0})\right)^{2}},
		\end{eqnarray*}
		and this finish the proof. Therefore, we assume  for any fixed positive constant $a_{0}\in (|\cos\theta|, \frac{3n+1}{5n-1})$, it holds
		\begin{eqnarray*}
			\frac{|\nabla u|^{2}}{v}-\cos\theta\sin\beta u_{\beta}\geq (1-a_{0})u_{1}, 
		\end{eqnarray*}
		then it follows
		\begin{eqnarray*}
			I_{13}+I_{23}&=&-\left[\frac{n(1+Kd)|\nabla u|^{4}}{2v^{2}}(\rho+\frac{1}{\rho})-\frac{n(1+Kd)|\nabla u|^{2}}{2v}(\rho+\frac{1}{\rho})\cos\theta\sin\beta u_{\beta}\right]\\
			&&-\left[\frac{n\cos\theta \sigma(\nabla u,\nabla d)}{2}
			(\rho+\frac{1}{\rho})(\frac{|\nabla u|^{2}}{v}-\cos\theta\sin\beta \nabla_{\partial_{\beta}}u)\right]\\
			&=&-\frac{n}{2}u_{1}S(\rho+\frac{1}{\rho})(\frac{|\nabla u|^{2}}{v}-\cos\theta \sin\beta u_{\beta})\\
			&\leq& -\frac{n}{2}(1-a_{0})Su_{1}^{2}(\rho+\frac{1}{\rho}).
		\end{eqnarray*}
		Since $|\cos\theta|<a_{0}$ and we choose $\varepsilon:=\frac{\varepsilon_{0}}{2}\in (0,1)$ with $\varepsilon_{0}:=\frac{3n+1-a_{0}(5n-1)}{4n(1-a_{0})}>0$, we have $(n-1)(1+a_{0})-4(1-\varepsilon)(1-a_{0})n<0$, then
		\begin{eqnarray*}
			&&I_{13}+I_{23}+I_{12}+I_{22}+I_{32}\\
			&\leq& -\frac{n}{2}(1-a_{0})Su_{1}^{2}(\rho+\frac{1}{\rho})+\frac{1-r_{0}^{2}}{2r_{0}}\frac{(n-1)(1+|\cos\theta|)S}{16(1-\varepsilon)}(\rho-\frac{1}{\rho})^{2}\rho e^{U}u_{1}^{2}\\
			&\leq& u_{1}^{2}S\left[\frac{(n-1)(1+|\cos\theta|)}{16(1-\varepsilon)}
			(\rho-\frac{1}{\rho})^{2}\frac{2\rho}{\rho^{2}+1}-\frac{n}{2}(1-a_{0})(\rho+\frac{1}{\rho})\right]\\
			&=&\frac{u_{1}^{2}S}{8\rho(\rho^{2}+1)(1-\varepsilon)}\Big[\left((n-1)(1+a_{0})-4n(1-\varepsilon)(1-a_{0})\right)(\rho^{4}+1)\\
			&&-\left(2(n-1)(1-a_{0})+8n(1-\varepsilon)(1-a_{0})\right)\rho^{2}
			\Big]\\
			&\leq& -\alpha_{0}u_{1}^{2},
		\end{eqnarray*}
		where $\alpha_{0}$ is a positive constant, which depends on $n, a_{0}, ||u||_{C^{0}}$. 
		
		Adding all above terms into \eqref{3.4}, we deduce
		\begin{eqnarray*}
			0&\leq& -\frac{\varepsilon_{0}}{2}\frac{2r_{0}}{1-r_{0}^{2}}\frac{(1+Kd)}{\rho v^{2}e^{U}}\sum\limits_{\alpha=2}^{n}u_{\alpha\alpha}^{2}-\alpha_{0}u_{1}^{2}+O(\frac{1}{v})\sum\limits_{\alpha=2}^{n}|u_{\alpha\alpha}|+O(v)\\&\leq &-\alpha_0 u_1^2+O(v),
		\end{eqnarray*}
		which follows  
		\begin{eqnarray*}
			u_{1}\leq C. 
		\end{eqnarray*}
		where the positive constant $C$ depends only on $n, r_{0}$, and $\|u\|_{C^{0}}$. Hence we complete the proof.
	\end{proof}
    Following the above argument for the case $\HH^{n+1}$,  we can get the uniform height and gradient estimates for the scalar parabolic equation \eqref{scalar equ in hyperbolic} in the case $\SS^{n+1}$. That is.
    \begin{proposition}\label{pro sphere gradient}
    If $u: \bar{\SS}^{n}_{+}\times [0, T)\rightarrow \RR$ solves \eqref{scalar equ in hyperbolic} and $G$ has the form \eqref{G sphere},  $|\cos\theta|<\frac{3n+1}{5n-1}$, then 
    	\begin{eqnarray}
    		\|u\|_{C^{1}\big(\bar{\SS}_{+}^{n}\times[0,T)\big)}\leq C,
    	\end{eqnarray} 
    where the constant $C$ is a positive constant, depending on the intial datum.
    \end{proposition} 
	For the concise of this paper, we leave the proof of Proposition \ref{pro sphere gradient} to the interested readers.
	
	In conclusion, we have the following convergence for the flow \eqref{flow with capillary in hyperbolic} both in hyperbolic space and spherical space.
	\begin{proposition}\label{3.4}
		The smooth solution of flow \eqref{flow with capillary in hyperbolic}  exists for all time and has uniform $C^\infty$-estimates, if the initial hypersurface $\Sigma_0\subset \bar B_R\subset M^{n+1}(K)$ with $K=\pm 1$ is star-shaped in the sense of Definition \ref{star-shaped}  and   $|\cos\theta|<\frac{3n+1}{5n-1}$.
	\end{proposition}
	\begin{proof}
		Proposition \ref{pro3.2},   \ref{pro3.3} and \ref{3.4}  say that $u$ is uniformly bounded in $C^1(\bar\SS^n_+\times[0,T))$, then the scalar equation in \eqref{scalar equ in hyperbolic} is uniformly parabolic.
		 Since $|\cos\theta|<1$, hence the desired conclusion follows from the standard quasi-linear parabolic theory with strictly oblique boundary condition theory (cf. \cite{Gary1996,LSU}).
	\end{proof}	
	Finally, we obtain  the convergence result by using the argument  in \cite{SWX,WeX21}, that is,  we complete the proof of Theorem \ref{thm1.1}.
	\begin{proposition}\label{pro3.4}
		If the initial hypersurface $\S_0\subset  \bar B_R\subset M^{n+1}(K)$ is star-shaped capillary boundary hypersurface and $|\cos\theta|<\frac{3n+1}{5n-1}$, then the flow \eqref{flow with capillary in hyperbolic}  smoothly converges to a uniquely determined spherical cap $C_{\theta, r}(a)$ given by \eqref{static model} with capillary boundary, as $t\rightarrow +\infty$.
	\end{proposition}
	\begin{proof}
		In the following,  we present a complete proof of the convergence for flow \eqref{flow with capillary in hyperbolic} in the case $\HH^{n+1}$. Since the proof for the case  $\SS^{n+1}$ is similar, we omit it here.
		
		From the proof of Corollary \ref{cor1.3} and uniform $C^\infty$-estimate, we see
		\begin{eqnarray}
			\int_{0}^{\infty}\int_{\Sigma_{t}}\sum\limits_{1\leq i<j\leq n}(\kappa_{i}(x,t)-\kappa_{j}(x,t))^{2}\bar{g}(X_{a},\nu)dA_{t} \leq C,\notag	
		\end{eqnarray}
		where the $\kappa_{i}(x,t), i=1,\cdots, n$ are the principal curvatures of the radial graph at the point $(x,t)$.  Together with the uniform estimate, we see $\bar{g}(X_{a},\nu)$ and $dA_{t}$ are uniformly bounded, it follows
		\begin{eqnarray*}
			\max\limits_{\substack{x\in \Sigma_{t} \\ 1\leq i< j\leq n}}|\kappa_{i}(x,t)-\kappa_{j}(x,t)|=o_{t}(1),
		\end{eqnarray*}
		where $o_{t}(1)$ denotes a quantity which goes to zero as $t\to +\infty$. Hence any convergent subsequence of $x(\cdot, t)$ converges to a spherical cap as $t\to +\infty$.
		
		Next, we show that the limit spherical cap is unique  by  following the argument in \cite{SWX,WeX21}. First, we know  any convergent subsequence of $x(\cdot, t)$   smoothly converges to a spherical cap $C_{\theta ,\rho_{\infty}}(a_{\infty})$. Since the volume is preserved along with the flow \eqref{flow with capillary in hyperbolic},  the radius is independent of the choice of the subsequence of $t$. Now we just need to show that $a_{\infty}=a$.
		
		Denote $\rho(\cdot, t)$ be the radius of the unique spherical 
		cap centered at the point  $\sqrt{\rho^{2}(\cdot, t)+r_{0}^{2}+2\rho(\cdot, t)r_{0}\cos\theta}a$ with contact angle $\theta$ passing through the point $x(\cdot, t)$. Following from the same barrier argument in Proposition \ref{pro3.2}, 
		\begin{eqnarray}
			\rho_{\max}(t):=\max\limits_{x\in M}\rho(x, t)=\rho(\xi_{t}, t),\notag
		\end{eqnarray}  
		is non-increasing with respect to $t$, for some point $\xi_t\in M$, hence the limit $\lim\limits_{t\rightarrow +\infty}\rho_{\max}(t)$ exists and it is clear that $\rho_{\max}(t)\geq \rho_{\infty}$.  We claim that
		\begin{eqnarray}\label{3.5}
			\lim\limits_{t\rightarrow +\infty}\rho_{\max}(t)=\rho_{\infty}.
		\end{eqnarray}
		We prove the above claim by a contradiction. Suppose \eqref{3.5} is not true, then there exists a constant $\varepsilon>0$, when $t$ is large enough, such that
		\begin{eqnarray}\label{3.6}
			\rho_{\max}(t)>\rho_{\infty} +\varepsilon.
		\end{eqnarray}
		By the definition of $\rho(\cdot, t)$,
		\begin{eqnarray}\label{3.7}
			2\<x, a\>\sqrt{\rho^{2}+r_{0}^{2}+2\rho r_{0}\cos\theta}=|x|^{2}+r_{0}^{2}+2\rho r_{0}\cos\theta,
		\end{eqnarray}
		taking the time derivative on the both sides for \eqref{3.7}, we get
		\begin{eqnarray*}
			\langle x_{t}, x-\sqrt{\rho^{2}+r_{0}^{2}+2\rho r_{0}\cos\theta}a\rangle=\left(\frac{(\rho+r_{0}\cos\theta)\<x, a\>}{\sqrt{\rho^{2}+r_{0}^{2}+2\rho r_{0}\cos\theta}}-r_{0}\cos\theta\right)\p_t {\rho}.
		\end{eqnarray*}
		We evaluate at  point $(\xi_{t}, t)$, note that $\Sigma_{t}$ is tangential to $C_{\theta,\rho_{\max}}(a)$ at $(\xi_{t}, t)$, it implies 
		\begin{eqnarray}
			(\nu_{0})_{\Sigma_{t}}(\xi_{t}, t)=(\nu_{0})_{\partial C_{\theta,\rho_{\max}(a)}}=\frac{x-\sqrt{\rho_{\max}^{2}(t)+r_{0}^{2}+2\rho_{\max}(t)r_{0}\cos\theta}a}{\rho_{\max}(t)},\notag
		\end{eqnarray}
		therefore
		\begin{eqnarray}
			\begin{aligned}\label{3.8}
				&\left(\frac{(\rho_{\max}(t)+r_{0}\cos\theta)\<x, a\>}{\sqrt{\rho_{\max}^{2}(t)+r_{0}^{2}+2\rho_{\max}(t) r_{0}\cos\theta}}-r_{0}\cos\theta\right)\partial_{t}\rho_{\max}(t)\\
				&=e^{-U}F\left\<(\nu_{0})\big|_{\Sigma_{t}}(\xi_{t}, t), x-\sqrt{\rho_{\max}^{2}(t)+r_{0}^{2}+2\rho_{\max}(t)r_{0}\cos\theta}a\right\>\\
				&=e^{-U}\rho_{\max}(t) \left(nV_{a}+n\sinh R\cos \theta \bar{g}(Y_{a},\nu)-H\bar{g}(X_{a}, \nu)\right).
			\end{aligned}	
		\end{eqnarray}
		Since the spherical $C_{\theta, \rho_{\max}}(a)$ is a static solution to flow \eqref{flow with capillary in hyperbolic},  the mean curvature $\bar{H}$ of $C_{\theta,\rho_{\max}}(a)$ in $(\bar{B}^{\HH}_{R}, \bar{g})$ is 
		\begin{eqnarray*}
			\bar{H}=e^{-u}\left[\frac{n}{\rho_{\max}(t)}+\frac{2n}{1-|x|^{2}}\<x, \nu_{0}\big|_{\partial C_{\theta,\rho_{\max}}(a)}\>\right]
			=\frac{n\left(1-r_{0}^{2}-2r_{0}\rho_{\max}(t)\cos\theta\right)}{2\rho_{\max}(t)}, 
		\end{eqnarray*}
		then  
		\begin{eqnarray}\label{3.9}
			\frac{nV_{a}+n\sinh R\cos\theta {g}(Y_{a},\nu)}{\bar{g}(X_{a},\nu)}\Big|_{C_{\theta,\rho_{\max}} (a)}=\frac{n\left(1-r_{0}^{2}-2r_{0}\rho_{\max}(t)\cos\theta\right)}{2\rho_{\max}(t)}.
		\end{eqnarray}
		Since	$x(\cdot, t)$ converges to $C_{\theta, \rho_{\infty}}(a_{\infty})$ and $\rho_{\infty}$ is uniquely determined, we have
		\begin{eqnarray}\label{3.10}
			\bar	H\rightarrow \frac{n(1-r_{0}^{2}-2\rho_{\infty}r_{0}\cos\theta)}{2\rho_{\infty}}.
		\end{eqnarray}
		We claim that there exist a positive constant $\delta >0$, such that
		\begin{eqnarray}\label{3.11}
			\frac{\left(\rho_{\max}(t)+r_{0}\cos\theta\right)\<x,a\>}{\sqrt{\rho_{\max}^{2}(t)+r_{0}^{2}+2\rho_{\max}(t) r_{0}\cos\theta}}-r_{0}\cos\theta\geq \delta.
		\end{eqnarray}
		In fact, by \eqref{3.7}, we have
		\begin{eqnarray}\label{3.12}
			\begin{aligned}
				&\<x, a\>^{2}(\rho^{2}+2\rho r_{0}\cos\theta+r_{0}^{2})\\
				&=\frac{1}{4}(|x|^{2}+r_{0}^{2})^{2}+\rho r_{0}\cos\theta(|x|^{2}+r_{0}^{2})+\rho^{2}r_{0}^{2}\cos^{2}\theta,
			\end{aligned}
		\end{eqnarray}
		combining with  \eqref{3.7}, it yields 
		\begin{eqnarray*}
			&&(\rho+r_{0}\cos\theta)\<x,a\>-r_{0}\cos\theta\sqrt{\rho^{2}+2\rho r_{0}\cos\theta+r_{0}^{2}}\\
			&=&(\rho+r_{0}\cos\theta)\<x, a\>-r_{0}\cos\theta\frac{|x|^{2}+r_{0}^{2}+2\rho r_{0}\cos\theta}{2\<x, a\>}\\
			&=&\frac{1}{\rho \<x,a\>}\left[\rho(\rho+r_{0}\cos\theta)\<x,a\>^{2}-\rho^{2}r_{0}^{2}\cos^{2}\theta-\rho r_{0}\cos\theta\frac{|x|^{2}+r_{0}^{2}}{2}
			\right]\\	&=&\frac{1}{\rho\<x, a\>}\left[\frac{(|x-r_{0}a||x+r_{0}a|)^{2}}{4}+\frac{1}{2}\rho r_{0}\cos\theta(|x|^{2}+r_{0}^{2})-\rho r_{0}\cos\theta \<x, a\>^{2}
			\right],
		\end{eqnarray*}
		together with  Proposition \ref{pro3.2}, it yields that Claim \eqref{3.11} is true.
		
		On the other hand, by \eqref{3.5}, \eqref{3.9}, \eqref{3.10} and the uniform estimates we established before, then there exists some large constant $T_{0}$ satisfying for $t>T_{0}$, it holds
		\begin{eqnarray*}
			e^{-U}\left(nV_{a}+n\sinh R\cos\theta\bar{g}(Y_{a}, \nu)-H\bar{g}(X_{a}, \nu)\right)\big|_{x(\xi_{t}, t)}\leq -C\varepsilon.
		\end{eqnarray*}
		Finally, by \eqref{3.8} and \eqref{3.11}, we conclude that there exists a positive constant $C_{0}$  such that
		\begin{eqnarray*}
			\frac{d}{dt}\left(\rho_{\max}(t)\right)\leq -C_{0}\varepsilon.
		\end{eqnarray*}
		This contradicts to the fact that $\lim\limits_{t\rightarrow+\infty}\frac{d}{dt}\left(\rho_{\max}(t)\right)=0$, so \eqref{3.5} is true. Similarly, one can obtain
		\begin{eqnarray}
			\lim\limits_{t\rightarrow+\infty}\rho_{\min}(t)=\rho_{\infty}.\notag
		\end{eqnarray} 
		Therefore,  $\lim\limits_{t\rightarrow+\infty}\rho(\cdot, t)=\rho_{\infty}$. This implies any limit of the convergent subsequence is the spherical cap $C_{\theta, \rho_{\infty}}(a)$ around $a$ with radius $\rho_\infty$. We complete the proof of Proposition \ref{pro3.4}, which follows also Theorem \ref{thm1.1}.
	\end{proof}

	\textbf{Acknowledgment.}
	Both authors would like to express sincere gratitude to Prof. Xinan Ma and Prof. Guofang Wang for their constant
	encouragement and many inspiring conversations in this subject.
	XM is partially supported by   CSC (No. 202106340053) and	the doctoral dissertation creation project of USTC.  LW is partially supported by China 
	Postdoctoral Science Foundation (No. 2021M702143) and NSFC (No. 12201003, 12171260).

\end{document}